# RANDOM GROWTH MODELS WITH POLYGONAL SHAPES[1]

By Janko Gravner and David Griffeath

*University of California, Davis and University of Wisconsin*

We consider discrete-time random perturbations of monotone cellular automata (CA) in two dimensions. Under general conditions, we prove the existence of half-space velocities, and then establish the validity of the Wulff construction for asymptotic shapes arising from finite initial seeds. Such a shape converges to the polygonal invariant shape of the corresponding deterministic model as the perturbation decreases. In many cases, exact stability is observed. That is, for small perturbations, the shapes of the deterministic and random processes agree exactly. We give a complete characterization of such cases, and show that they are prevalent among threshold growth CA with box neighborhood. We also design a nontrivial family of CA in which the shape is exactly computable for all values of its probability parameter.

**1. Introduction.** Discrete local models for random growth and deposition have been a staple of rigorous research in probability since the Hammersley and Welsh paper [18] on first passage percolation about 40 years ago. Apart from their role as a testing ground for probabilistic techniques, a voluminous physics literature [26, 28] testifies to their importance in understanding the evolution of natural systems far from equilibrium. The most basic tool, introduced in [18] and ubiquitous ever since, is *subadditivity*: the process dominates one restarted from an already occupied point. Clearly, this imposes a monotonicity property on the model, but, as we will see, not much more. The result is the existence of an *asymptotic shape*: started from a finite seed, and scaled by time, the occupied set converges to a deterministic convex limit. Elegant as this method is, it is nonconstructive and as a result fails to provide any detailed information about the limiting set. Thus asymptotic properties of subadditive sequences are still an active area of

Received August 2003; revised April 2005.

[1]Supported in part by NSF Grants DMS-02-04376 and DMS-02-04018.

research [2, 34]. Are there cases when the shape can be exactly identified? Research on this topic has so far primarily focused on growth from *infinite* initial states, also known as *random interfaces*. Methods have ranged from hydrodynamic limits based on explicit identification of invariant measures [30, 31], to techniques arising from exactly solvable systems in mathematical physics [14, 20], and to perturbation arguments based on supercritical oriented percolation [6, 8] which imply that some interfaces move with the speed of their deterministic counterparts. For other related rigorous and empirical results see [16, 22, 27].

The main aim of this paper is to extend the perturbation approach to show that the finite limit shape of a random growth model may also agree with that of a deterministic one. At issue is not merely whether small random errors induce small changes (we will see that this is always the case), but rather whether the shape can stay exactly the same. This property, which we call *exact stability*, is only valid under substantial assumptions, as the model has to have opposite structure, in an appropriate sense, from the additive one considered in [6]. In the process, we extend the result of [3] to obtain the Wulff characterization of the invariant shape. We also show that exact stability is far from rare; in fact, almost all members of arguably the most natural family of two-dimensional growth models, the threshold growth cellular automata with square neighborhood, are exactly stable. Finally, we show how to employ exactly solvable systems to construct *one* example which has a computable shape for every value of its probability parameter. (Although they are invaluable in suggesting universal phenomena, exactly solvable examples are extremely difficult to come by.)

The random rules we describe below can be thought of as discrete counterparts to the KPZ equation, which in turn is touted as a universal scaling model for any local growth and deposition process in physics [26, 28], in particular for crystal growth [28]. This we mention because the well-studied *roughening transition* in crystallography, whereby a crystal loses its polygonal shape as the ambient temperature increases, produces pictures which are strikingly similar to ours [32]. While this transition is usually thought to be an equilibrium phenomenon, the present results at least suggest that it may have a *dynamic* counterpart.

We now proceed to precise formulations. Unfortunately, these require a large number of definitions related to our previous work. Although we do not use any of the results from [10] explicitly, a glance at that paper's first two sections may help to motivate what follows.

Our basic framework consists of two-state cellular automata (CA). In general, such a cellular automaton is specified by the following two ingredients. The first is a finite *neighborhood* $\mathcal{N} \subset \mathbf{Z}^2$ of the origin, its translate $x + \mathcal{N}$ then being the neighborhood of point $x$. By convention, we assume that $\mathcal{N}$ contains the origin. Typically, $\mathcal{N} = B_\nu(0, \rho) = \{x : \|x\|_\nu \leq \rho\}$, where $\|\cdot\|_\nu$ is the $\ell^\nu$-norm. When $\nu = 1$ the resulting $\mathcal{N}$ is called the range $\rho$ Diamond neighborhood, while if $\nu = \infty$ it is referred to as the range $\rho$ Box neighborhood. (In particular, range 1 Diamond and Box neighborhoods are also known as von Neumann and Moore neighborhoods, resp.) The second ingredient is a map $\pi : 2^{\mathcal{N}} \to \{0, 1\}$, which flags the *sufficient configurations* for occupancy. More precisely, for a set $A \subset \mathbf{Z}^2$, we define $\mathcal{T}(A) \subset \mathbf{Z}^2$ by adjoining every $x \in \mathbf{Z}^2$ for which $\pi((A_t - x) \cap \mathcal{N}) = 1$. Then, for a given *initial subset* $A_0 \subset \mathbf{Z}^2$ of occupied points, we define $A_1, A_2, \ldots$ recursively by $A_{t+1} = \mathcal{T}(A_t)$. Accordingly, occupied and vacant sites will often be denoted by 1's and 0's, respectively. Our main focus will be starting states $A_0$ which consist of a possibly large, but finite set of 1's surrounded by 0's. However, we will also consider other initial states, namely half-spaces and wedges.

We restrict to two-dimensional dynamics for two main reasons. First, almost every step in higher dimensions introduces new technical complications, some quite serious. In fact, there are new phenomena, and the classification of Theorem 2 below becomes much more complex. Second, some of our techniques are intrinsically two-dimensional, such as the explicitly solvable example of Section 6, the lattice geometry and analytic number theory of Section 7, and even combinatorial properties studied in [3]. Nevertheless, some results—notably Theorem 1—do readily generalize to arbitrary dimension.

Our key assumption is that the CA dynamics are *monotone* (or *attractive*), that is, $S_1 \subset S_2$ implies $\pi(S_1) \leq \pi(S_2)$. Note that specifying a monotone dynamics is the same as specifying an *antichain* of subsets of $\mathcal{N}$: the inclusion minimal sets $S$ with $\pi(S) = 1$ having the property that none of them is a subset of another. Surprisingly, the number of possible monotone dynamics (known as a Dedekind number) is possible to estimate for large $\mathcal{N}$. Some typical properties of monotone CA are also known [23]. Unfortunately, it turns out that for large box neighborhoods the asymptotic proportion of supercritical rules (see the definition below) is negligible. Other interesting properties seem to present great difficulties. In studying typical monotone CA rules, it is therefore desirable to restrict to a simpler class.

A natural such class consists of *totalistic* monotone CA, those for which $\pi(S)$ depends only on the cardinality $|S|$ of $S$. In other words, there exists a *threshold* $\theta \geq 0$ such that $\pi(S) = 0$ whenever $|S| < \theta$ and $\pi(S) = 1$ whenever $|S| \geq \theta$. This much studied case is also known by the name *threshold growth* (TG) CA.

Induced by $\mathcal{T}$ is a growth transformation $\bar{\mathcal{T}}$ on closed subsets of $\mathbf{R}^2$, given by

$$\bar{\mathcal{T}}(B) = \{x \in \mathbf{R}^2 : 0 \in \mathcal{T}((B - x) \cap \mathbf{Z}^2)\}.$$

In words, one translates the lattice so that $x \in \mathbf{R}^2$ is at the origin, and applies $\mathcal{T}$ to the intersection of Euclidean set $B$ with the translated lattice.

It is easy to verify that the two transformations are *conjugate*,

$$\mathcal{T}(B \cap \mathbf{Z}^2) = \bar{\mathcal{T}}(B) \cap \mathbf{Z}^2.$$

It will become immediately apparent why $\bar{\mathcal{T}}$ is convenient. Let $S^1 \subset \mathbf{R}^2$ be the set of unit vectors and let

$$H_u^- = \{x \in \mathbf{R}^2 : \langle x, u \rangle \leq 0\}$$

be the closed half-space with outward normal $u \in S^1$. Then there exists a $w(u) \in \mathbf{R}$ so that

$$\bar{\mathcal{T}}(H_u^-) = H_u^- + w(u) \cdot u$$

and consequently

$$\mathcal{T}^t(H_u^- \cap \mathbf{Z}^2) = (H_u^- + tw(u) \cdot u) \cap \mathbf{Z}^2.$$

If $w(u) > 0$ for every $u$ we call the CA *supercritical*. A supercritical CA hence enlarges every half-space. This is equivalent to existence of a finite set $A_0$ which *fills space*, that is, $\bigcup_{t \geq 0} A_t = \mathbf{Z}^2$ [3, 10]. *All initial sets will be assumed to fill space from now on.* Set

$$K_{1/w} = \bigcup \{[0, 1/w(u)] \cdot u : u \in S^1\}$$

and let $L$ be the polar transform of $K_{1/w}$, that is,

$$L = K_{1/w}^* = \{x \in \mathbf{R}^2 : \langle x, u \rangle \leq w(u)\}.$$

Then one can prove the following limiting shape result for any finite $A_0$:

$$\lim_{t \to \infty} \frac{A_t}{t} = L,$$

where the limit is taken in the Hausdorff metric. In short, the shape $L = L(\pi)$ is obtained as the Wulff transform of the speed function $w : S^1 \to \mathbf{R}$, which for small neighborhoods is readily computable by hand or by computer. Furthermore, $L$ is always a polygon and the Hausdorff distance between $A_t$ and $tL$ is bounded in time $t$ [9, 10, 11, 12, 13, 37].

To formulate the stability properties of $L$ under random perturbations, we begin by introducing a general monotone random dynamics. The function $\pi$ differs from the one described above in that it has values in $[0, 1]$. Upon seeing a set of occupied sites $x + S$ in its neighborhood at time $t$, a site becomes occupied at time $t + 1$ independently with probability $\pi(S)$. To obtain a monotone rule we require that $\pi(S_1) \leq \pi(S_2)$ whenever $S_1 \subset S_2$.

More precisely, introduce i.i.d. vectors $\xi_{x,t}$, $x \in \mathbf{Z}^2$, $t = 0, 1, 2, \ldots$, with $2^{|\mathcal{N}|}$ coordinates $\xi_{x,t}(S)$, which are Bernoulli($\pi(S)$) for every $S \subset \mathcal{N}$. We assume that these are coupled so that $\xi_{x,t}(S_1) = 1$ implies that $\xi_{x,t}(S_2) = 1$

whenever $S_1 \subset S_2$. The construction of such a coupling is left as an exercise for the reader. The random sets $A_1, A_2, \ldots$ are now determined by

$$A_{t+1} = \{x : \xi_{x,t}(((x+\mathcal{N}) \cap A_t) - x) = 1\}.$$

To avoid some trivialities and inessential complications, we assume that 1's only grow *by contact*: $\pi(\varnothing) = 0$, and that $\pi$ is *symmetric*: $-\mathcal{N} = \mathcal{N}$ and $\pi(-S) = \pi(S)$. Much more substantial is the assumption that $\pi$ *solidifies*: $\pi(S) = 1$ whenever $0 \in S$. These three properties, together with monotonicity, will be our *standing assumptions* throughout the paper.

For every random $\pi$, we set $p = \min\{\pi(S) : \pi(S) > 0\}$, define the associated deterministic dynamics by its map $\pi_d(S) = \mathbb{1}_{\{\pi(S)>0\}}$, and label the iteration transform $\mathcal{T}$ as before. We will say that $\pi$ is a *p-perturbation* of the CA $\mathcal{T}$. For many purposes the *standard p-perturbation*, which has $\pi(S) = p$ whenever $\pi(S) > 0$, suffices.

We say that a *p*-perturbation of $\mathcal{T}$ has *shape* $L_\pi$ if

$$\lim_{t \to \infty} \frac{A_t}{t} = L_\pi$$

almost surely, in the Hausdorff metric, for every finite initial set $A_0$ which fills space. We say that $\mathcal{T}$ has *exactly stable* shape $L$ if there exists a $p < 1$ such that $L_\pi = L$ (which of course subsumes the existence of the shape $L_\pi$) for the standard, and hence any, *p*-perturbation $\pi$. For a standard perturbation, we also write $L_p = L_\pi$. Thus $L_1 = L$. Recall that the deterministic growth at time $t$ is included in a constant fattening of $tL_1$; hence the same is true of any *p*-perturbation.

As already mentioned, such considerations are in the general direction of the vintage Durrett–Liggett flat edge result [6]. To describe their result in our context, recall that a deterministic CA is *additive* if $\pi(S)$ equals 1 precisely when $S$ is nonempty. In this case $K_{1/w} = \mathcal{N}^*$ and $L = \text{co}(\mathcal{N})$. Moreover, any standard perturbation is a first passage percolation model, and as such has an almost sure (deterministic) limiting shape $L_p$ for each $p > 0$ [5, 29]. For the von Neumann neighborhood, Durrett and Liggett proved that, if $p$ is close to 1, then $L_p$ is close to $L$ and in fact inherits from $L$ flat edges in the four diagonal directions. However, they show that $L_p$ is *not* equal to $L$, due to the fact that its extent in the coordinate directions is strictly less than 1.

The existence of a limiting shape $L_\pi$ for general random dynamics does not immediately follow from standard subadditivity arguments. A sufficient condition is a property of $\mathcal{T}$ we call *local regularity*. Namely, for every initial state $A_0$ there exists a constant $C$ so that the following is true for every fixed (deterministic) assignment of $\xi_{x,t}$: every $x \in A_t$ at distance at least $C$ from $A_0$ has an occupied set $G \subset A_t$ entirely within distance $C$ of $x$ such that $G$ fills space.

Note that local regularity is a combinatorial condition involving every possible way $A_t$ can evolve, and thus has nothing to do with probability. At first it seems a condition not likely to be often satisfied, but the opposite appears true. One can easily check local regularity directly for many cases with small $\mathcal{N}$, and it holds generally for box neighborhood TG CA. All known counterexamples involve "strange" neighborhoods [3]. Under this condition, it can readily be shown that $L_\pi$ exists.

Besides finite shapes, limiting profiles from half-spaces are of considerable interest. The first reason is that their Monte Carlo approximations can be computed much more efficiently (see Remark 2 in Section 8). The second is that they are important for shapes from other infinite sets, such as wedges and holes [13]. For finite seeds also, the Wulff transform (see Corollary 1.1 below), which expresses the asymptotic shape in terms of half-space velocities, is very handy. However, the limit theorem in [3] does not extend to infinite seeds, as restarting requires an a priori upper bound on fluctuations. Here we provide the missing step, which establishes the following large deviations bound, referred to as the *Kesten property* in [13]. (See [21] for a similar result in the first passage context.)

THEOREM 1. *Let $\pi$ be a p-perturbation of a locally regular supercritical CA and let the initial set be $A_0 = H_u^- \cap \mathbf{Z}^2$ for $u \in S^1$. Then there exists a deterministic $w_\pi(u) > 0$ such that*

$$H_u^- + t(w_\pi(u) - \varepsilon) \cdot u \subset A_t \subset H_u^- + t(w_\pi(u) + \varepsilon) \cdot u$$

*within the lattice ball of radius $t^2$ with probability at least $1 - \exp(-c_\varepsilon t)$. Here $c_\varepsilon > 0$ as soon as $\varepsilon > 0$.*

COROLLARY 1.1. *For a p-perturbation $\pi$ of a locally regular CA and finite initial sets which fill space,*

$$\frac{A_t}{t} \to L_p = K^*_{1/w_\pi},$$

*in the Hausdorff metric, almost surely.*

Next is a generalization of the flat edge result [6]. In particular this implies that $L_p \to L_1$ when $p \to 1$, as promised.

PROPOSITION 1.2. *Given a standard p-perturbation of a locally regular CA and any $\varepsilon > 0$, there exists a $p < 1$ close enough to 1 that $L_p$ agrees with $L_1$ outside the $\varepsilon$-neighborhood of the set of corners of $L_1$.*

Our second theorem provides necessary and sufficient conditions for exact stability. Before its statement, it is instructive to look at the three supercritical Moore TG CA. The $\theta = 1$ case is additive and exact stability cannot hold. (This can be proved by the methods of [6] or [25], but we give a different argument in Section 3.) For $\theta = 2$ one finds that $K_{1/w} = \text{co}(\mathcal{N})$ and hence this is a *quasi-additive* case, that is, a CA with convex $K_{1/w}$. Quasi-additive CA share many properties with additive ones [10, 11, 13], and lack of exact stability turns out to be among them. Finally, in the $\theta = 3$ case $K_{1/w}$ has 16 vertices, of which three successive ones are $(0,1)$, $(1,2)$, $(1,1)$, and the remaining 13 are then continued by symmetry. (This set, which the reader is invited to compute, is the innermost region of Figure 7.) Eight of these are the only points that $K_{1/w}$ shares with the boundary of its convex hull. In a sense, the fact that these eight vertices form a discrete set makes this CA as unlike a quasi-additive one as possible. This turns out to be precisely the condition needed for exact stability.

Accordingly, we denote

$$\partial K' = \partial(K_{1/w}) \cap \partial(\text{co}(K_{1/w})),$$

and describe the relevance of properties of this set in our main result.

THEOREM 2. *Consider a supercritical locally regular CA (which also satisfies our standing assumptions) given by $\mathcal{T}$, with limiting shape $L_1$, and its standard p-perturbation. There are three possibilities:*

Case 1. *$\partial K'$ consists of isolated points, no three of which are collinear.*
   *Then the following hold for $p < 1$ close enough to 1:*
   (S1) *$L_p = L_1$.*
   (S2) *Convergence to $L_1$ is tight: for any $\varepsilon > 0$, there exists an $M$ so that, for any $t$ and $x \in tL_1$, $P(x$ is within $M$ of $A_t) \geq 1 - \varepsilon$.*
   (S3) *There exists a large $C$ so that, with probability 1, $(t - C \log t)L_1 \cap \mathbf{Z}^2 \subset A_t$ eventually.*

Case 2. *$\partial K'$ consists of isolated points, three of which are collinear.*
   *Then (S1) and (S3) still hold for $p < 1$ close enough to 1, but tightness (S2) no longer does. Instead, for any $p < 1$ there exists a $c > 0$ so that a corner of $tL_1$ is eventually at distance at least $c \log t$ from $A_t$, a.s.*

Case 3. *$\partial K'$ includes a line segment.*
   *Then (S1) no longer holds. Instead, for any $p < 1$ there exists a $c > 0$ so that a corner of $tL_1$ is at distance at least $ct$ from $A_t$, a.s.*

Figure 1 shows a box neighborhood TG example for each of the three cases, from left to right, with periodic shading of updates: Case 1 (range 1, $\theta = 3$, $p = 0.9$), Case 2 (range 2, $\theta = 7$, $p = 0.95$) and Case 3 (range 2, $\theta = 8$, $p = 0.95$).

The fundamental difference between Moore $\theta = 2$ and $\theta = 3$ TG CA is their mistake-fixing ability, which we now illustrate. Suppose we start each case with a large copy of the invariant shape and remove a finite chunk of occupied sites at the boundary. Regardless of the location of such a hole, the $\theta = 3$ case eventually repairs (or "erodes") it and thus the hole's effect is bounded in time. Figure 2 provides a demonstration. This eroding property can be used to favorably compare the random dynamics on infinite wedges, determined by the corners of $L_1$, to *Toom rules* [35]. The corners are then patched together by an oriented percolation comparison in the middles of the edges. In Case 2, mistakes are still fixed, but for wider wedges than in Case 1, and corners must be rounded off accordingly.

By contrast, the $\theta = 2$ TG CA can only repair holes away from the corners, while those at the corners have a lasting effect, as also seen in Figure 2. In a random dynamics, such mistakes pile up and induce a linear slowdown.

Given the exact stability criterion of Theorem 2, it is natural to ask whether a typical supercritical CA has an exactly stable shape or not. As already mentioned, properties of typical monotone CA seem difficult to characterize. We will thus restrict our attention to a special family, TG CA with range $\rho$ box neighborhoods $\mathcal{N}_\rho$. These are supercritical for $\theta \leq \rho(2\rho+1)$ [10]. The smallest examples are already illuminative. As $\mathcal{N}_1$ has already been discussed, $\mathcal{N}_2$ is next in line and turns out to have $\theta = 1, 2, 3, 5, 8$ in Case 3, $\theta = 7, 9, 10$ in Case 2, and $\theta = 4, 6$ in Case 1. For *very* large ranges, $\theta$'s in Case 3 form a small minority, as the following theorem demonstrates.

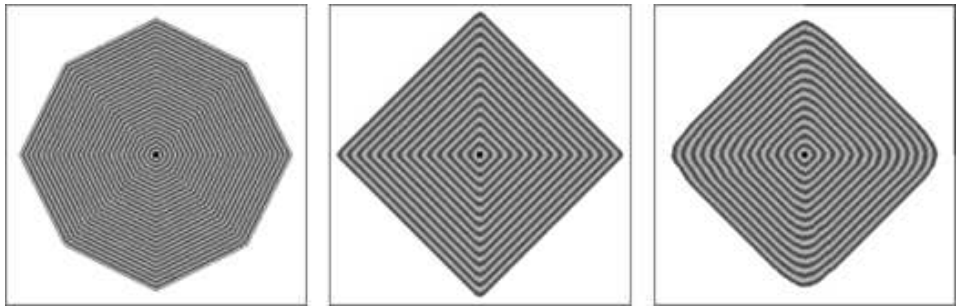

FIG. 1. *The three cases of Theorem 2.*

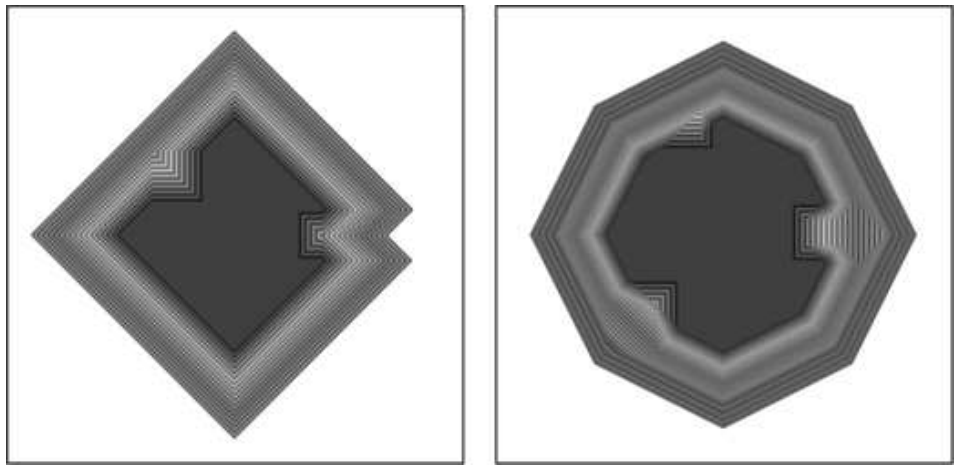

FIG. 2. *Error correcting for $\theta = 2$ and $\theta = 3$.*

THEOREM 3. *Fix an arbitrary $\varepsilon > 0$. Among all supercritical range $\rho$ box neighborhood TG CA, the proportion of those which are not exactly stable is for large $\rho$ between $1/\log^{h+\varepsilon} \rho$ and $1/\log^h \rho$. Here $h = 2(1 - 1/\log 2 - \log \log 2 / \log 2) \approx 0.172$.*

The proof of Theorem 3 connects the number of $\theta$'s which lack exact stability to the number of distinct products of pairs of natural numbers between 1 and $\rho$. This latter is known as the Linnik–Vinogradov–Erdős problem, for which sharp asymptotic bounds were given by Hall and Tenenbaum [17]. We have no result on the division between Cases 1 and 2, but conjecture that Case 2 is much more prevalent.

The rest of the paper is organized as follows. Section 2 contains the proof of a slightly weaker version of Theorem 1 and its Corollary 1.1. Section 3 deals with Case 3, while Section 4 lays the geometric groundwork for the remaining cases and proves Proposition 1.3. In Section 5 we introduce Toom's method and complete the proof of Theorem 2. Section 6 is devoted to a single example for which we can compute the shape for *all* values of the probability parameter $p$. In Section 7 we take a closer look at the collection of $K_{1/w}$'s for fixed-range box neighborhoods, an analysis which culminates with the proof of Theorem 3. Finally, in Section 8 we finish the proof of Theorem 1 and discuss other related issues in lesser detail.

**2. Proof of Theorem 1.** Recall that Theorem 1 deals with supercritical locally regular CA and their $p$-perturbations. These will be our context throughout this section. We will allow all constants $C$ and $c$ to depend on $\mathcal{T}$ and $p$ in addition to their explicitly stated dependencies. (We emphasize that these constants will not, however, depend on the direction $u$.) In this section we only obtain a lower bound of the form $1 - \exp(-c_\varepsilon t / \log^2 t)$ on the probability of the event in Theorem 1.

Many times below we will *restart* the random dynamics at a deterministic time or a random stopping time $\tau$. This simply means that only $\xi_{x,t}$ with $t \geq \tau$ are used, with an initial state at time $\tau$ which will be specified.

LEMMA 2.1. *Assume that a finite $A_0 \ni 0$ fills the plane. Assume that $x$ is at distance $n$ from the origin. Then there exist constants $c, C > 0$ (depending on $A_0$) so that $P(x \notin \mathcal{T}^k(A_0)) \leq e^{-ck}$ for $k \geq Cn$.*

PROOF. Call $x$ *surrounded* at time $t$ if $x + A_0 \subset A_t$. By supercriticality, there exists a time $t_0$ at which $\pm e_1$ and $\pm e_2$ are all surrounded in the deterministic dynamics. Let $C_0 = |\mathcal{T}^{t_0}(A_0)|$. If $p_0 = p^{C_0}$, then $\pm e_1$ and $\pm e_2$ are all surrounded at time $t_0$ with probability at least $p_0$. Take a shortest lattice path $\wp : 0 = x_0, x_1, \ldots, x_n = x$. We now define i.i.d. geometric$(p_0)$ random variables $T_1, \ldots, T_n$ as follows. Run the dynamics for time $t_0$. If $x_1$ is surrounded at this time, $T_1 = 1$, otherwise restart the dynamics with $A_0$ at time $t_0$. Now run the restarted dynamics for time $t_0$; if it surrounds $x_1$ at this time, $T_1 = 2$, otherwise restart again with $A_0$, and so on. In general, on the event $\{T_i = k\}$, $T_{i+1}$ is the minimal $\ell \geq 1$ for which the dynamics restarted at time $k + (\ell - 1)t_0$ with $x_i + A_0$ surrounds $x_{i+1}$ at time $t_0$.

By monotonicity and exponential Chebyshev,

$$P(x \notin \mathcal{T}^{kt_0}(A_0)) \leq P(T_1 + \cdots + T_n \geq k) \leq e^{-\lambda k} E(\exp(\lambda T_1))^n,$$

for any $\lambda > 0$. To conclude the proof, choose $\lambda$ small enough that $E(\exp(\lambda \times T_1)) < \infty$. □

Note that this lemma implies that $w_\pi(u)$, if it indeed exists, is bounded away from 0 uniformly in $u$, for any $p > 0$.

LEMMA 2.2. *Assume that $|\mathcal{T}(A_0) \setminus A_0| = n$, and start the p-perturbation from the same initial set $A_0$. If $\tau = \inf\{t : \mathcal{T}(A_0) \subset A_t\}$, then $E(\tau) \leq p^{-1} \times (\log n + 3)$.*

PROOF. Note that all such sites attempt to get occupied simultaneously, each of them at each time with probability $p$. Hence $\tau$ is geometric$(p)$ for $n = 1$. For $n \geq 2$, write $a = -\log(1-p)$ and divide the sum below into terms with $k \leq a^{-1} \log n$ and with $k > a^{-1} \log n$ to obtain

$$E(\tau) = \sum_{k=0}^{\infty} (1 - (1 - e^{-ak})^n)$$

$$\leq a^{-1} \log n + 1 + \sum_{i=0}^{\infty} (1 - (1 - e^{-ai} n^{-1})^n)$$

$$\leq a^{-1} \log n + 1 - \sum_{i=0}^{\infty} n \cdot \log(1 - e^{-ai} n^{-1})$$

$$\leq a^{-1} \log n + 1 + 2 \sum_{i=0}^{\infty} e^{-ai}$$

$$= a^{-1} \log n + 1 + 2(1 - e^{-a})^{-1},$$

and $p = 1 - e^{-a} < a$. □

PROOF OF THEOREM 1 WITH WEAKER PROBABILITY ESTIMATE. Without loss of generality, we can assume that $u$ lies on or above $y = |x|$, that is, $\langle u, e_2 \rangle \geq 1/\sqrt{2}$. Let $\mathcal{F}_t = \sigma\{\xi_{x,s} : s \leq t-1, x \in \mathbf{Z}^2\}$, $t = 1, 2, \ldots$.

Let $T_n$ be the first time $(0, n)$ becomes occupied started from $H_u^-$ and set $\bar{T}_n = T_n \wedge Cn$. By Lemma 2.1, $P(T_n \neq \bar{T}_n) \leq e^{-cn}$, for a large enough $C$ and some $c > 0$. The crucial step is this $L^\infty$ bound:

$$(2.1) \qquad |E(\bar{T}_n | \mathcal{F}_{s+1}) - E(\bar{T}_n | \mathcal{F}_s)| \leq C' \log n,$$

for any $s \leq Cn$ and some constant $C'$.

Recall that $\bar{T}_n$ is a deterministic function of $\xi_{x,t}$ where $(x,t)$ ranges over all space-time sites. As $\mathcal{N}$ is finite, $\bar{T}_n$ depends only on a small subset of these variables. To be more precise, let $\mathcal{L}_n$ comprise the sites $(x,t)$ for which $\bar{T}_n$ depends on $\xi_{x,t}$. Then $|\mathcal{L}_n| \leq Cn^3$ and we can assume that the filtration ignores all other sites. At time $s \leq Cn$, let $\partial A_s$ consist of all the sites outside $A_s$ which would become occupied if the deterministic dynamics were applied to $A_s$. Trivially, $|\partial A_s| \leq |\mathcal{L}_n|$.

Restart the dynamics at time $s+1$ with $A_s$. Let $\tau_s$ be the waiting time after this at which *all* sites in $\partial A_s$ are occupied, that is, $\tau_s = \inf\{k : \partial A_s \subset A_{s+1+k}\}$. By Lemma 2.2, $E(\tau_s | \mathcal{F}_s) \leq C'' \log n$.

We now prove (2.1). We will repeatedly use the strong Markov property and monotonicity of the dynamics. To get the lower bound in (2.1), assume the worst case: no sites outside $A_s$ (i.e., in $\partial A_s$) get occupied, and therefore the dynamics faces an unchanged situation at time $s+1$. Therefore,

$$E(\bar{T}_n | \mathcal{F}_{s+1}) \leq E(\bar{T}_n | \mathcal{F}_s) + 1.$$

For the upper bound, assume that $\mathcal{F}_{s+1}$ reveals that all sites in $\partial A_s$ get occupied. Before we know $\mathcal{F}_{s+1}$, we can only assume this happens after time $\tau_s$, and so the dynamics with the additional information is dominated by the one restarted at time $s + \tau_s$ from the occupied set $A_s \cup \partial A_s$. It follows that

$$E(\bar{T}_n | \mathcal{F}_s) \leq E(\bar{T}_n | \mathcal{F}_{s+1}) + E(\tau_s | \mathcal{F}_s) \leq E(\bar{T}_n | \mathcal{F}_{s+1}) + C'' \log n.$$

This proves (2.1).

Now let $a_n = E(T_n)$, $\bar{a}_n = E(\bar{T}_n)$. By (2.1) and Azuma's inequality [19, 33],

$$(2.2) \qquad P(|\bar{T}_n - \bar{a}_n| > s) \leq 2\exp(-cs^2/(n \log^2 n)).$$

However,

$$|a_n - \bar{a}_n| \leq E(T_n \mathbb{1}_{\{T_n \geq Cn\}}),$$

which is bounded by Lemma 2.1. From this it follows that

$$P(|T_n - a_n| > s) \leq P(|\bar{T}_n - \bar{a}_n| > s/2 - C) + P(T_n - \bar{T}_n > s/2),$$

and after another application of Lemma 2.1 and suitable redefinition of $c$,

(2.3) $\qquad P(|T_n - a_n| > s) \leq 2\exp(-cs^2/(n\log^2 n)) + e^{-cs}.$

For an integer $i$, let $y_i$ be the largest $j$ for which $(i,j) \in H_u^-$. Then let $T_n'$ be the first time at which *all* sites in $B' = \{(i, y_i + n) : |i| \leq n^2\}$ are occupied. Moreover, let $T_n''$ be the first time at which all the sites $B'' = \{(i,j) : |i| \leq n^2, y_i + n - C \leq j \leq y_i + n\}$ are occupied, where $C$ is the diameter of the neighborhood $\mathcal{N}$. Restart the dynamics at time $T_n'$ with the occupied set at this time. Note that local regularity implies that within a constant time the deterministic dynamics occupies a large ball within a constant distance of any occupied point. By monotonicity, the deterministic dynamics would occupy $B''$ in $t_1$ additional time steps, where $t_1$ is a constant which only depends on $\mathcal{T}$. Applying Lemma 2.2 $t_1$ times one thus obtains

$$E(T_n'' - T_n') \leq C \log n.$$

Furthermore, let $T_n'(i)$ be the time the dynamics reaches $(i, n + y_i)$ and let $T_n(i)$ be the first time $(i, n + y_i)$ becomes occupied *from the modified initial set* $y_i e_2 + H_u^-$. The reason for this convoluted condition is that $T_n(i)$ with the same $n$ are identically distributed, but this is not true for $T_n'(i)$.

To deal with different starting sets for $T_n'(i)$, let $S_n$ be the time the random dynamics fills $H_u^- \cap B(0, n^2)$ from $-e_2 + H_u^-$ (which is contained in all starting sets). By a similar argument as in the previous paragraph

$$E(S_n) \leq C \log n.$$

Therefore, with $a_n'(i) = E(T_n'(i))$,

$$0 \leq a_n'(i) - a_n \leq E(S_n) \leq C \log n.$$

Furthermore, the argument leading to (2.3) can be carried out with $T_n$ replaced by $T_n'(i)$ and $a_n$ by $a_n'(i)$.

Therefore, for $s > C \log n$,

$$P(T_n - T_n' \geq s) \leq P(|T_n - a_n| > s/4) + \sum_{|i| \leq n^2} P(|T_n'(i) - a_n(i)| > s/4)$$

$$\leq Cn^2 \exp(-cs^2/n \log^2 n) + Cn^2 e^{-cs}.$$

It follows that

$$E(T_n - T_n') \leq C\sqrt{n}\log^2 n + \int_{C\sqrt{n}\log^2 n}^{\infty} P(T_n - T_n' \geq s)\, ds,$$

which after a short computation implies that $E(T_n - T_n') \leq C\sqrt{n}\log^2 n$.

We are almost done, but need an estimate for yet another approximation to $T_n$. Let $T_n'''$ be the first occupation time of $(0,n)$ started from $B'' - ne_2 = \{(i,j) : |i| \leq n^2, y_i - C \leq j \leq y_i\}$. Then, for $0 < k < Cn$,

$$P(T_n''' - T_n > k) \leq P(T_n''' \neq T_n) \leq P(T_n > Cn) \leq e^{-cn},$$

while for $k \geq Cn$,

$$P(T_n''' - T_n > k) \leq P(T_n''' \geq k) \leq e^{-ck}$$

by Lemma 2.1. Hence $E(T_n''' - T_n)$ is bounded above by a constant.

Now assume that $0 \leq m \leq n$. Restarting the growth process at time $T_n$, we get

$$a_{m+n} \leq a_m + a_n + E(T_n'' - T_n) + E(T_m''' - T_m) \leq a_m + a_n + C\sqrt{n}\log^2 n.$$

By the deBruijn–Erdős subadditive theorem [33], $a_n/n$ converges to a finite positive number $a$, which of course depends on $p$ and $u$. We declare $w_\pi(u) = \langle u, e_2 \rangle / a$.

To finish the proof, take first an $(i,j)$ outside $H_u^- + tw_\pi(u)(1+\varepsilon) \cdot u$. Let $n = j - y_i \geq tw_\pi(u)(1+\varepsilon)/\langle u, e_2 \rangle = t(1+\varepsilon)/a$. Then

$$\begin{aligned}
P((i,j) \in A_t) &= P(T_n'(i) \leq t) \\
&\leq P(T_n'(i) \leq na/(1+\varepsilon)) \\
&\leq P(|T_n'(i) - a_n'(i)| \geq na\varepsilon/2) \\
&\leq \exp(-cn/\log^2 n),
\end{aligned}$$

for a large enough $n$. This proves the weaker version of the upper bound in Theorem 1. The lower bound is proved similarly. $\square$

Several remarks are in order. First, note that the proof avoids the subadditive ergodic theorem altogether, by combining properties of subadditive sequences with large deviation estimates.

Second, it is in fact possible, by the same methods, to obtain a superadditive relation for $a_n$ of the same order, namely,

$$a_{m+n} \geq a_m + a_n - C\sqrt{n}\log^2 n.$$

A closer look at the proof of the deBruijn–Erdős theorem (from [33]) then gives a rate of convergence for $a_n$: $|a_n - a| = O(\log^2 n/\sqrt{n})$, which can be used to show that, within a lattice ball of radius $t^2$, $A_t$ is a.s. between $(t \pm C\sqrt{t}\log^3 t) \cdot w_\pi(u) \cdot u + H_u^-$.

Third, the proof uses supercriticality and regularity only to "fill in." For any monotone, local, interface solidification with automatic coherence the proof remains valid. While we will not attempt to precisely define the concept, automatic coherence certainly holds when the interface moves upward (i.e., $u = e_2$) and the growth is such that an empty site can never have an occupied site directly above it. Perhaps the simplest example is the random dynamics in which a site becomes occupied for sure with two or more occupied neighbors in its von Neumann neighborhood and with probability $p$ with an occupied site directly below. Another class of examples are the $K$-exclusion processes [31]. For some of these examples, the fluctuation estimates mentioned above may be new.

Finally, and curiously, there seems no way to make the proof work for general monotone dynamics which do not solidify. Such cases thus remain an intriguing challenge.

LEMMA 2.3. *Fix an $a > 0$ and $\varepsilon > 0$. Then there exist constants $c, C$ so that the following holds. Start the dynamics from $A_0$ consisting of sites inside $(H_u^- \setminus (-Cu + H_u^-)) \cap B(0, Cn)$. Then, $A_n$ includes all sites inside $B = ((H_u^- + nw_\pi(u)(1 - \varepsilon)u) \setminus H_u^-) \cap B(0, Cn)$ with probability at least $1 - e^{-cn/\log^2 n}$.*

PROOF. Let $T'(x)$ [resp. $T(x)$] be the first occupation time of $x \in B$ started from the stated $A_0$ (resp. from $H_u^-$). By Lemma 2.1, $P(\sup_{x \in B} T(x) > Cn) \leq e^{-cn}$. However, by a "speed of light" argument, on $\{\sup_{x \in B} T(x) \leq Cn\}$ the equality $T(x) = T'(x)$ holds for all $x \in B$. The claim now follows from Theorem 1. □

LEMMA 2.4. *The function $w_\pi : S^1 \to \mathbf{R}$ is continuous.*

PROOF. Again assume that $\langle u, e_2 \rangle \geq 1/\sqrt{2}$. For a fixed large $C$ and small $\varepsilon > 0$, $H_v^- - Ct\varepsilon e_2 \subset H_u^- \subset H_v^- + Ct\varepsilon e_2$, within the lattice ball of radius $Ct$, provided $\|u - v\|_2 < \varepsilon/2$. Let $E_v(k, t)$ be the event that all sites on the $y$-axis up to $k$ are occupied at time $t$ started from $H_v^-$.

By Lemma 2.3 and Theorem 1, both the events $E_v((1 - C\varepsilon)w_\pi(u)t/\langle v, e_2\rangle, t)$ and $E_v((1 + \varepsilon)w_\pi(v)t/\langle v, e_2\rangle, t)^c$ happen with probability (very) close to 1. This is only possible if $(1 - C\varepsilon)w_\pi(u) \leq (1 + \varepsilon)w_\pi(v)$. An analogous reverse inequality is proved similarly. □

PROOF OF COROLLARY 1.1. An $\varepsilon > 0$ will be fixed throughout this proof.

For any direction $u$, Theorem 1 implies that with probability exponentially close to 1

$$A_t \subset H_u^- + tw_\pi(u)(1 + \varepsilon) \cdot u.$$

It follows that with probability 1

$$\frac{A_t}{t} \subset H_u^- + w_\pi(u)(1 + \varepsilon) \cdot u,$$

eventually. This is therefore true simultaneously for any finite collection of $u$'s and then by Lemma 2.4,

$$\frac{A_t}{t} \subset \bigcap_{u \in S^1} H_u^- + w_\pi(u)(1 + \varepsilon)^2 \cdot u = (1 + \varepsilon)^2 K_{1/w_\pi}^*$$

eventually.

For the lower bound, take a bounded, strictly convex, $\mathcal{C}^2$ set $K_\varepsilon \supset K_{1/w_\pi}(1+\varepsilon)$. Then $L_\varepsilon = K_\varepsilon^*$ is $\mathcal{C}^2$ and has for small enough $\delta > 0$ the property described in the next paragraph.

Start with $A_0$ consisting of sites inside $nL_\varepsilon$. Take $k = n^2/\delta$ Euclidean points $x_0, \ldots, x_{k-1}$ on the boundary of $nL_\varepsilon$, chosen so their directions are equidistant vectors in $S^1$, and let $u_0, \ldots, u_{k-1}$ be the outside normals to $nL_\varepsilon$ at the chosen points. The enlarged set

$$L_n(\delta) = \bigcap_{i=0}^{k-1} x_i + (1-\delta)\sqrt{n} w_\pi(u_i) u_i + H_{u_i}^-$$

includes $(n + \sqrt{n})L_\varepsilon$.

Now run the random dynamics from $A_0$ for $\sqrt{n}$ time steps. Since $L_\varepsilon$ has $\mathcal{C}^2$ boundary, we need to go just a constant distance inside to "see" the relevant portion of $H_{u_i}^-$. To be more precise, $((-Cu_i + H_{u_i}^-) \setminus (-2Cu_i + H_{u_i}^-)) \cap B(x_i, C\sqrt{n})$ is included in $nL_\varepsilon$, for all $i$. By Lemma 2.3, with probability at least $1 - \exp(-c\sqrt{n}/\log^2 n)$, all the sites in $(n + \sqrt{n})L_\varepsilon$ become occupied.

Repeat the above procedure (running the random dynamics for $\sqrt{n}$ time steps) $3\sqrt{n}$ times. As a result, $(n + j\sqrt{n})L_\varepsilon \subset A_{n+j\sqrt{n}}$ for $j = 1, \ldots, 3\sqrt{n}$ (in particular, $4nL_\varepsilon \subset A_{4n}$), with probability at least $p_n = 1 - \exp(-c\sqrt{n}/\log^2 n)$.

Now fix an $a < 1$ and find a large $k_0$ so that $\prod_{n \geq 2^{2k_0}} p_n > a$. Let $T_0$ be the first time $2^{2k_0} L_\varepsilon \subset A_t$. By what we proved so far,

$$P((2^{2k} + j2^k)L_\varepsilon \subset A_{2^{2k}+j2^k - T_0} \text{ for } j = 0, \ldots, 3 \cdot 2^k, k = 1, 2, \ldots) > a.$$

We thus have a strictly increasing sequence of integers $b_m$ with $b_{m+1} - b_m = o(b_m)$, such that

$$b_m(1 - \varepsilon) L_\varepsilon \subset A_{b_m}$$

eventually, with probability at least $a$, thus a.s., as $a$ was arbitrary. For any $t$ between $b_m$ and $b_{m+1}$,

$$t(1 - \varepsilon)^2 L_\varepsilon \subset b_m(1 - \varepsilon) L_\varepsilon \subset A_{b_m} \subset A_t$$

eventually, finishing the proof of the lower bound. $\square$

## 3. Lack of exact stability in Case 3.

Fix a $u \in S^1$. Let $\ell_u$ be the boundary line of $-w(u) \cdot u + H_u^-$. Note that $w(u)$ is the largest number $h > 0$ for which $\pi((-h \cdot u + H_u^-) \cap \mathcal{N}) = 1$. Therefore, $\mathcal{N} \cap \ell_u$ must contain at least one site.

In general, for any line $\ell$ in the plane which does not go through the origin, let its *open* (resp. *closed*) *lower cut* $\mathcal{L}^o(\ell)$ [resp. $\mathcal{L}^-(\ell)$] be the set of points in $\mathcal{N}$ which lie in the open (resp. closed) half-space of $\ell^c$ which does not contain the origin. We emphasize here (as this convention will be used extensively) that the points in $\mathcal{L}^o(\ell)$ will be called *below* the line $\ell$, and that *left* and *right* on the line are from the perspective of an observer who stands on $\ell$ and looks toward the origin.

We will make good use of duality between lines in $K_{1/w}$ and points of $\mathcal{N}$ in the sequel. The next lemma is our first example of this duality. To illustrate its statement (as well as the introduced terminology), let us consider an example. Assume that we are dealing with a TG CA and fix a direction $u$. Suppose also that $\ell_u$ contains an $x_u \in \mathcal{N}$ such that a line $\ell$ obtained by a small rotation of $\ell_u$ around $x_u$ has exactly $\theta - 1$ sites in $\mathcal{L}^o(\ell)$. Note that for a sufficiently small rotation no other site but $x_u$ is in $\ell \cap \mathcal{N}$. (An example for the range 2 Box TG CA with $\theta = 8$ is depicted on the right-hand side of Figure 3.) Therefore, for $v$ close enough to $u$, $\ell_v$ is obtained by rotation of $\ell_u$ around $x_u$. A little geometric argument involving polar coordinates then shows that the boundary of $K_{1/w}$ must be flat at $u/w(u)$. The lemma makes a stronger and more general statement and is illustrated by the left-hand side of Figure 3.

LEMMA 3.1. *The following are equivalent for a $u \in S^1$.*

(1) *There exists a line through $u/w(u)$ which in a small neighborhood of $u/w(u)$ lies in $K_{1/w}$.*
(2) *There exists a point $x_u \in \ell_u \cap \mathcal{N}$ so that if $\ell$ is a line through $x_u$ and is a rotation of $\ell_u$ by a small enough angle, $\pi(\mathcal{L}^o(\ell)) = 0$.*

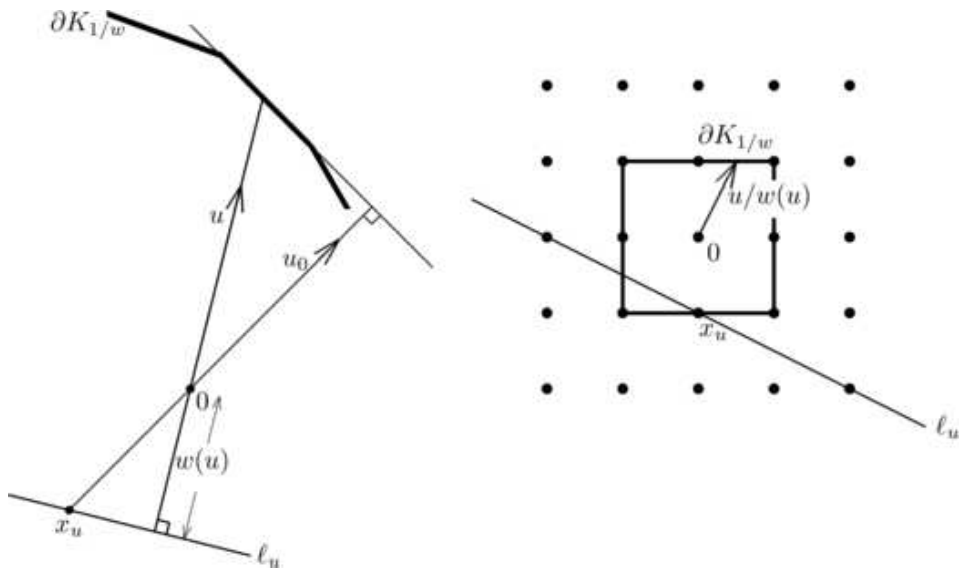

FIG. 3. *Illustration of Lemma 3.1 and its proof.*

In case $\partial K_{1/w}$ is locally a line at $u/w(u)$, $x_u$ in (2) is unique. In fact, the smaller angle between $\partial K_{1/w}$ and $u/w(u)$ is the same as the smaller angle between the vector $x_u$ and $\ell_u$.

PROOF. Note that a short line segment through $u/w(u)$ perpendicular to the vector $u_0 \in S^1$ is given in polar coordinates (with the angle represented by a unit vector $v$) by the collection of vectors

$$\left\{ \frac{\langle u, u_0 \rangle}{w(u)} \cdot \frac{v}{\langle v, u_0 \rangle} : \|v - u\| < \alpha \right\},$$

for a small $\alpha > 0$.

Assume first that the statement (2) holds. Let $u_0 = -x_u/\|x_u\|$ be the unit vector pointing from $x_u$ to the origin. Then (2) says that for $v$ close enough to $u$, $w(v) \leq w(u)\langle v, u_0 \rangle / \langle u, u_0 \rangle$. It follows that

$$(3.1) \qquad \frac{1}{w(v)} \geq \frac{1}{w(u)} \frac{\langle u, u_0 \rangle}{\langle v, u_0 \rangle}.$$

The polar representation of a line mentioned above immediately demonstrates the implication $(2) \Rightarrow (1)$.

To prove the reverse implication, note that (1) implies that (3.1) holds for some $u_0$, and let $x'_u = -w(u)u_0/\langle u, u_0 \rangle$. Then $x'_u$ has the properties required of $x_u$, except it may not lie in $\mathcal{N}$. However, we can let $x_u$ be the closest site in $\ell_u \cap \mathcal{N}$. (We can in fact go in either direction from $x_u$.) The fact that $\mathcal{N}$ is discrete ensures that a parallel translation from $x'_u$ to $x_u$ of any line $\ell$ close to $\ell_u$ does not pass through any site of $\mathcal{N}$. Thus (2) is satisfied.

To prove the last statement, note that two different $x_u$ would, by (3.1), produce two distinct open line segments, which would meet at $u/w(u)$ and which would both be included $K_{1/w}$. But then a flat portion of $\partial K_{1/w}$ near $u/w(u)$ would be impossible. $\square$

LEMMA 3.2. *Fix a $u \in S^1$ which satisfies the condition of Lemma 3.1 and pick a corresponding $x_u$. For $v$ close enough to $u$, the concave wedge $Q = H_u^- \cup H_v^-$ satisfies $\bar{\mathcal{T}}(Q) \subset -x_u + Q$. When $\partial K_{1/w}$ is locally a line at $u/w(u)$, $Q$ is invariant: $\bar{\mathcal{T}}(Q) = -x_u + Q$.*

PROOF. The first part follows from Lemma 3.1: a point in $\bar{\mathcal{T}}(Q) \cap (-x_u + Q)^c$ would imply that a point on the boundary of $-x_u + Q$ sees a sufficient configuration in the interior of $Q$, but clearly $-x_u$ is in the most advantageous position for this. This would translate, for $\ell$ as in (2) of Lemma 2.1, into $\pi(\mathcal{L}^o(\ell) \cup \mathcal{L}^o(\ell_u)) = 1$, but if the rotation is sufficiently small (by discreteness of $\mathcal{N}$) $(\mathcal{L}^o(\ell) \cup \mathcal{L}^o(\ell_u)) \cap \mathcal{N} = \mathcal{L}^o(\ell) \cap \mathcal{N}$, a contradiction.

The second part also follows because in this case $\pi(\mathcal{L}^-(\ell)) = 1$. For, otherwise $x_u$ could be moved to the next point in $\ell_u \cap \mathcal{N}$ for which $\pi(\mathcal{L}^-(\ell)) = 1$.

[Again, such a point must exist or else $w(u)$ could be decreased.] This would contradict uniqueness. □

When $u$ satisfies the assumption of Lemma 3.2 there exists an invariant wedge of the following form:
$$Q' = (-Mv_1 + H_{v_1}^-) \cup (-Mv_2 + H_{v_2}^-) \cup H_u^-,$$
where $v_1$ and $v_2$ are close to, but on different sides of, $u$ and $M$ is large enough. In particular, a hole of shape $Q'$ dug into $H_u^-$ may be translated by the dynamics, but is never filled. If the creation of such holes is random they pile up and, as we will demonstrate by the comparison process we now introduce, slow down the interface.

The following randomly growing surface will be useful here and in Section 8. At every time $t = 0, 1, 2, \ldots$ a site $x \in \mathbf{Z}$ has a height $\eta_t'(x) \in \mathbf{Z}_+$, with $\eta_0' \equiv 0$. We will use two versions, which we call *fast* and *slow*, of the rule for increase in heights. Let $b(x, t)$ be Bernoulli random variables with $P(b(x, t) = 1) = p'$. The slow version evolves according to the following rule:
$$\eta_{t+1}'(x) = \begin{cases} \eta_t'(x) + 1, & \text{if } b(x, t) = 1 \text{ and} \\ & \eta_t'(y) \geq \eta_t'(x) \text{ for all } y \text{ with } |y - x| \leq 1, \\ \eta_t(x), & \text{otherwise,} \end{cases}$$
while the fast version updates as follows:
$$\eta_{t+1}'(x) = \begin{cases} \eta_t'(x) + 1, & \text{if } b(x, t) = 1 \text{ or} \\ & \eta_t'(y) > \eta_t'(x) \text{ for some } y \text{ with } |y - x| \leq 1, \\ \eta_t(x), & \text{otherwise.} \end{cases}$$
Note that the reverse dynamics, $\eta_t(x) = t - \eta_t'(x)$ changes the version and replaces $p'$ by $1 - p'$. We will assume that $b(x, t)$ are not necessarily independent, but have finite range dependence in space: if either $t_1 \neq t_2$ or $|x_1 - x_2| > r$, then $b(x_1, t_1)$ and $b(x_2, t_2)$ are independent.

LEMMA 3.3.

(1) *For the slow version: Given any $p' > 0$, there exist an $\alpha > 0$ and $c > 0$ so that*
$$P(\eta_t'(x) \leq \alpha t) \leq e^{-ct}.$$

(2) *For the slow version: Given any $\varepsilon > 0$, there exist a large enough $p'$ and a $c > 0$ so that*
$$P(\eta_t'(x) \leq (1 - \varepsilon)t) \leq e^{-ct}.$$

(3) *For the fast version: Given any $\varepsilon > 0$, there exist a small enough $p'$ and a $c > 0$ so that*
$$P(\eta_t'(x) \geq \varepsilon t) \leq e^{-ct}.$$

PROOF. The proof of (1) and (2) is a last passage percolation argument. By [24] we can in fact assume that the random variables $b(x,t)$ are independent. Once the neighborhood condition ($\eta'_t(y) \geq \eta'_t(x)$ for all $y$ with $|y - x| \leq 1$) is satisfied, a site $x$ has to wait a geometric($p'$) number of time steps before it increases. Accordingly, let $g(x,s)$ be i.i.d. geometric with success probability $p'$. By a simple inductive argument, it follows that the first time $T_n(x)$ when $\eta'_s(x) = n \geq 1$ equals

$$\max\left\{\sum_{i=0}^{n-1} g(x_i, i) : x_{n-1} = x, x_{i+1} \in \{x_i - 1, x_i, x_i + 1\} \text{ for } 0 \leq i < n - 1\right\}.$$

Hence

$$P(\eta'_s(x) \leq n) = P(T_n(x) > s) \leq 3^n P\left(\sum_{i=0}^{n-1} g(0, i) > s\right).$$

By an elementary large deviation computation, we get, for a fixed $p' > 0$ and a small enough $\alpha > 0$, $P(\eta'_s(x) \leq \alpha s) \leq \exp(-cs)$, which implies (1). Another large deviation computation gives $P(\eta'_s(x) \leq (1 - \varepsilon)s) \leq \exp(-cs)$ for a fixed $\varepsilon > 0$ and $p'$ close enough to 1. Finally, (3) follows from (2) by reversal. □

PROOF OF THEOREM 2 IN CASE 3. Let $u \in S^1$ be a direction of an interior point in a line segment of $\partial K'$. Then $u/w(u)$ belongs to the interior of a line segment of $\partial K_{1/w}$ (satisfying the condition of Lemma 3.1) and the corner of $L$ which corresponds to the edge of $\text{co}(K_{1/w})$ containing $u/w(u)$ moves, in the deterministic case, with speed $w(u)$ in direction $u$. The following claim will therefore finish the proof. Start the random dynamics from $A_0 = H_u^- \cap \mathbf{Z}^d$. Then, for some $\alpha > 0$,

$$(3.2) \qquad P(A_t \subset (1 - \alpha)w(u)tu + H_u^-) \geq 1 - e^{-ct}.$$

For simplicity, rotate the space so that $u = e_2$. Recall that $2M$ is the width of the bottom edge of $Q'$. Assign $\tilde{\eta}_1(i) = 0$ if $A_1 \cap (iMe_1 + Q') = \varnothing$, and $\tilde{\eta}_1(i) = 1$ otherwise. In general, let $\eta_t(i)$ be the smallest $k$ for which $A_t \cap (iMe_1 - kx_u + Q') = \varnothing$. It is clear that $\tilde{\eta}_t$ is for $M$ large enough dominated by $\eta_t = t - \eta'_t$ where $\eta'_t$ is the slow version from Lemma 3.3 and $p' = (1 - p)^k$, for some appropriately large $k$. The range of dependence $r$ depends on $M$ and angles between $v_1$ and $u$ and $v_2$ and $u$, but is clearly finite. Therefore, (3.2) follows from Lemma 3.3(1). □

**4. Flat edges of shapes for $p$ close to 1.** The setup is the same as in the previous section. In Lemma 4.1 below, the direction of a rotation of a line $\ell$ is determined by the direction of motion of the outward normal to the half-space in $\ell^c$ which does not contain the origin.

The left-hand side of Figure 4 depicts a general situation in the statement and proof of the lemma, while the right-hand side again presents a TG CA example. This time the range 2 Box case has $\theta = 7$ (see Figure 6). Note that if $\ell$ is a rotation of $\ell_u$ around $x_u^\ell$ by a small negative angle, then there are exactly $\theta - 1 = 6$ sites below $\ell$. The same is true for rotations around $x_u^r$ by a small positive angle. This translates to two line segments on the boundary of $K_{1/w}$ which meet at $u/w(u)$ at a convex angle, of $45°$ in this case.

LEMMA 4.1. *For a $u \in S^1$, assume that near $u/w(u)$ the boundary of $K_{1/w}$ consists of two lines at interior angle below $\pi$. Then there exist $x_u^\ell, x_u^r \in \ell_u \cap \mathcal{N}$ with the following properties. If $\ell$ is a small rotation of $\ell_u$ either through $x_u^\ell$ by a negative angle, or through $x_u^r$ by a positive angle, then $\pi(\mathcal{L}^-(\ell)) = 1$. In fact, the smaller angles between $\partial K_{1/w}$ and $u/w(u)$ are the same as the smaller angles between $x_u^\ell$ and $\ell_u$ and between $x_u^r$ and $\ell_u$, if $x_u^\ell$ and $x_u^r$ are chosen to be furthest apart.*

PROOF. The argument is very similar to the one for Lemma 3.1. The local equation for $\partial K_{1/w}$ to the right of $u/w(u)$ is given by $v \mapsto (v/w(u)) \cdot (\langle u, u_0 \rangle / \langle v, u_0 \rangle)$, for a suitably chosen $u_0$. Then

$$\frac{1}{w(v)} < \frac{1}{w(u)} \frac{\langle u, u_0 \rangle}{\langle v, u_0 \rangle}, \tag{4.1}$$

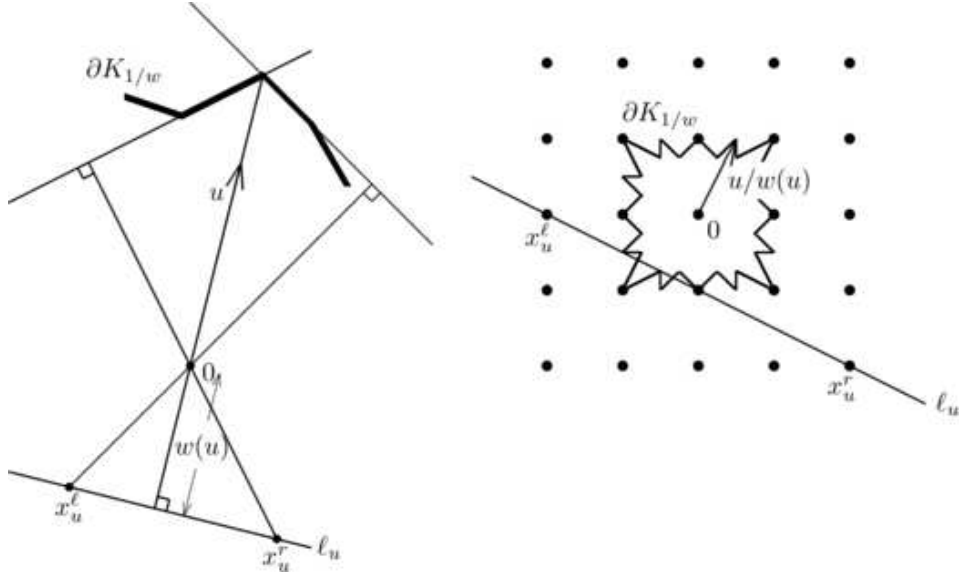

FIG. 4. *Illustration of Lemma 4.1.*

if $v$ is to the left of $u$. If $x'_u = -w(u)u_0/\langle u, u_0 \rangle$, this means that any small rotation $\ell$ of $\ell_u$ around $x'_u$ in the positive direction has $\pi(\mathcal{L}^o(\ell)) = 1$. Now simply move $x'_u$ rightward on $\ell_u$, to the first point on $\ell_u \cap \mathcal{N}$ for which $\pi(\mathcal{L}^o(\ell)) = 0$ for small positive rotations $\ell$. Such a point must exist, or else $\pi(\mathcal{L}^o(\ell_u)) = 1$, which contradicts the definition of $w(u)$. This defines $x^r_u$, which must be in $\mathcal{N} \cap \ell$, as small rotations contain no other sites of $\mathcal{N}$. The definition of $x^\ell_u$ is similar. □

LEMMA 4.2. *Fix $u \in S^1$, and $x^\ell_u, x^r_u$ as in Lemma 4.1, chosen as far apart as possible. If $v_1$ is a small positive rotation of $u$, then the convex wedge $W^+ = H^-_u \cap H^-_{v_1}$ is invariant: $\bar{\mathcal{T}}(W^+) = -x^r_u + W^+$. Similarly, if $v_2$ is a small negative rotation of $u$, then the convex wedge $W^- = H^-_u \cap H^-_{v_2}$ is invariant: $\bar{\mathcal{T}}(W^-) = -x^\ell_u + W^-$.*

PROOF. This proof is completely analogous to the one for Lemma 3.2. We omit the details. □

Note that the two wedges from Lemma 4.2 are moving toward each other. In particular, if we now dig any finite hole in $H^-_u$, it gets filled, as we state more precisely in the next corollary.

COROLLARY 4.3. *If $v_1$ and $v_2$ are as above, and*
$$A_0 = (H^-_u \cap (-Mv_1 + H^-_{v_1})) \cup (H^-_u \cap (-Mv_2 + H^-_{v_2})),$$
*then*
$$\mathcal{T}^t(A_0) = tw(u)u + H^-_u$$
*for $t \geq CM$.*

Assume now that $u_1, u_2 \in S^1$ are such that $u_1/w(u_1)$ and $u_2/w(u_2)$ are on the boundary of $\mathrm{co}(K_{1/w})$, and that $u_2$ is a positive (counterclockwise) rotation of $u_1$ (by the smaller angle between them). Assume also that for small positive (resp. negative) rotations $v$ of $u_2$ (resp. $u_1$), $v/w(v)$ is not on the boundary of $\mathrm{co}(K_{1/w})$. This assumption always holds in Cases 2 and 3 of Theorem 2 when $u_1/w(u_1)$ and $u_2/w(u_2)$ are vertices of $\mathrm{co}(K_{1/w})$, and also holds for other vectors in Case 2. Then
$$(w(u_1)u_1 + H^-_{u_1}) \cap (w(u_2)u_2 + H^-_{u_2}) = z + H^-_{u_1} \cap H^-_{u_2}$$
covers a vertex of $L$ [and when $u_1/w(u_1)$ and $u_2/w(u_2)$ are vertices of $\mathrm{co}(K_{1/w})$ also portions of corresponding edges of $L$]. The equation above defines the vector $z = z(u_1, u_2)$. While this wedge is by itself not necessarily invariant (although for small $\mathcal{N}$ it often is), a bounded perturbation with a suitably rounded corner is superinvariant [as in (3) of Lemma 4.5 below].

LEMMA 4.4. *Let $u_1$ and $u_2$ be as above. Apply Lemma 4.2 to $u = u_1$ to get the corresponding $W^-$ and $x_u^\ell$. The corner of $W^-$ moves slower than $H_{u_2}^-$, that is, $\langle -x_u^\ell, u_2 \rangle < w(u_2)$.*

PROOF. The conclusion is equivalent to $x_{u_1}^\ell \notin \mathcal{L}^-(\ell_{u_2})$. But this follows because $x_{u_1}^r \notin \mathcal{L}^o(\ell_{u'})$ for any $u'$ between $u_1$ and $u_2$. [If $K_{1/w}$ were a straight line between $u_1/w(u_1)$ and $u_2/w(u_2)$, $x_{u_1}^r$ would belong to $\ell_{u'}$ for any such $u'$.] □

An analogous version of Lemma 4.4 of course also holds for $u = u_2$ and the corresponding $W^+$.

LEMMA 4.5. *Assume $u_1$ and $u_2$ are as in Lemma 4.4. There exists a convex wedge $W_{u_1,u_2}$ which is included in, and outside a bounded neighborhood of the corner equal to, $H_{u_1}^- \cap H_{u_2}^-$, such that the following properties hold.*

(1) *When the part of the edge of $\mathrm{co}(K_{1/w})$ between $u_1/w(u_1)$ and $u_2/w(u_2)$ is completely included in $K_{1/w}$,*

$$\bar{\mathcal{T}}(W_{u_1,u_2}) = z + W_{u_1,u_2}.$$

(2) *When the open part of the edge of $\mathrm{co}(K_{1/w})$ between $u_1/w(u_1)$ and $u_2/w(u_2)$ has empty intersection with $K_{1/w}$,*

$$\bar{\mathcal{T}}(W_{u_1,u_2}) \supset z + ((W_{u_1,u_2} + B_2(0, \alpha)) \cap W_{u_1,u_2}),$$

*for some $\alpha > 0$.*

(3) *In every other case, $\bar{\mathcal{T}}(W_{u_1,u_2}) \supset z + W_{u_1,u_2}$.*

PROOF. The condition in (2) implies that every vector $v$ strictly between $u_1$ and $u_2$ has $x_{u_1}^r \in \mathcal{L}^o(\ell_v)$ and $x_{u_2}^\ell \in \mathcal{L}^o(\ell_v)$. This, together with Lemma 4.2, readily proves (2): one simply makes $W_{u_1,u_2}$ start with and end with wedges considered there, and connects them with a convex curve of small enough curvature.

The condition of (1) implies that every vector $v$ strictly between $u_1$ and $u_2$ has $x_{u_1}^r \in \ell_v$ and $x_{u_2}^\ell \in \ell_v$. Again, Lemma 4.2 implies that there are a succession of invariant wedges which connect the rays with normals $u_1$ and $u_2$. The final statement follows by a subdivision of the edge of $\mathrm{co}(K_{1/w})$ into subintervals of types considered in (1) or (2). □

COROLLARY 4.6. *Fix an $\varepsilon > 0$. There exists a convex set $L_\varepsilon \subset L$, which agrees with $L$ outside the $\varepsilon$-neighborhood of its corners, so that for large enough $M$,*

$$\bar{\mathcal{T}}(M \cdot L_\varepsilon) \supset (M+1)L_\varepsilon.$$

PROOF. The wedges from Lemma 4.5 can readily be combined to approximate an arbitrarily large multiple of $L$, within an error confined to a constant distance from its corners. □

Arguably, *oriented percolation* is the most useful comparison model in random spatial processes. We now introduce the version we will use. While this is in fact a random perturbation of a one-dimensional CA, it is, as the name suggests, best to think about it as a random occupied set which tries to establish long-range connections. Sites $(m,n) \in \mathbf{Z}_+ \times \mathbf{Z}_+$ are either occupied or empty ($n$ is often referred to as *a level*). The basic ingredients are Bernoulli($p'$) random variables $b(m,n)$, $m,n \geq 1$, such that $b(m_1,n_1)$ and $b(m_2,n_2)$ are independent whenever $|m_1 - n_1| > r$ or $n_1 \neq n_2$. (It is important that $r$ does not depend on $p'$.) Prescribe some occupied set in $\mathbf{Z}_+ \times \{0\}$. For $m \geq 0$ and $n \geq 1$, $(m,n)$ is occupied if $b(m,n) = 1$ and at least one of its *neighbors* $(m,n-1)$ and $(m-1,n-1)$ is occupied.

LEMMA 4.7. *Fix any $\alpha \in (0,1)$. Also fix a large integer $M$ and let $[0,M] \times \{0\}$ be occupied. If $p'$ is close enough to $1$, then for large enough $C = C(p')$ the probability of the following two events converges to $1$ as $M \to \infty$.*

(1) *Any $(m,n)$ with $\alpha n \leq m \leq (1-\alpha)n$ and $n = 0,1,\ldots$ is within distance $C \log n$ of an occupied point.*
(2) *For every $n$ there is a connection (through neighbors) of occupied points from level $n$ down to level $0$ which stays entirely in $\{(m,n) : n \geq 1, \alpha n \leq m \leq \alpha n + M + C \log n\}$.*

PROOF. These are standard applications of contour arguments (see, e.g., [5], so we omit the details. □

PROOF OF PROPOSITION 1.2. If $y_i, i = 0, \ldots, R-1$, are vectors pointing to the $R$ successive corners of $L$ in a counterclockwise order, then

$$(M+1)L_\varepsilon = \bigcup_{i=0}^{R-1} (y_i + L_\varepsilon).$$

We will now concentrate on the edge between $y_0$ and $y_1$. Start the random dynamics from $A_0 = M \cdot L_\varepsilon \cap \mathbf{Z}^2$. Say that $(m,n) \in \mathbf{Z}_+ \times \mathbf{Z}_+$ is *occupied* if all the lattice sites in $(n-m)y_0 + my_1 + ML_\varepsilon$ are at time $n$ included in $A_n$.

If $(m,n)$ is occupied, then both $(m,n+1)$ and $(m+1,n+1)$ are occupied with probability at least $p'$ which is a power of $p$ given by the number of lattice points in $(M+1)L_\varepsilon \setminus ML_\varepsilon$. Moreover, given any configuration of occupied points on level $n$, points on level $n+1$ which are $r$ apart are occupied independently. Here $r$ is any integer such that $(ry_0 + (M+1)L_\varepsilon) \cap (ry_1 + (M+1)L_\varepsilon) = \varnothing$.

The occupied points hence form an oriented percolation with a finite range of dependence. For $p'$ close enough to 1, Lemma 4.7(1) shows that the conclusion holds with probability which converges to 1 as $M \to \infty$. However, this is enough as every set is covered in almost surely finite time. □

**5. Exact stability in Cases 2 and 3.** We continue with the setup of the last two sections. Let us begin with a statement of Toom's theorem. We will only state the version we need here (for which the proof is in [35]), although the conclusion holds in considerably greater generality [4, 36]. A two-dimensional *Toom rule* is a deterministic CA $\mathcal{T}_T$ given by a map $\pi_T$ with the following property:

(T) There exists a line $\ell_T$ which does not go through the origin such that $\pi_T(S) = 1$ if and only if either $\mathcal{L}^o(\ell_T) \subset S$ or $\mathcal{N} \setminus \mathcal{L}^-(\ell_T) \subset S$.

Now introduce space-time *error sites*, those sites $(x,t)$ for which $b(x,t) = 0$. Here $b(x,t) \in \{0,1\}$, $x \in \mathbf{Z}^2$, $t = 1, 2 \ldots$, are assigned before the dynamics starts. The state of the *Toom rule with errors* is then given by $\eta_t \in \{0,1\}^{\mathbf{Z}^2}$, which satisfies $\eta_0 \equiv 1$ and

$$\eta_{t+1} = \mathcal{T}_T(\eta_t) \cap \{x : b(x, t+1) = 1\}, \qquad t = 0, 1, \ldots.$$

To develop some intuition, note that without errors a finite island of 0's in a sea of 1's gets eroded by $\mathcal{T}_T$, as it is "squeezed" between two half-spaces with boundaries parallel to $\ell_T$. (However, this island may move in the process.) Thus (T) is often called the *eroder condition*. A natural question is what happens with such a rule under persistent introduction of low-density errors.

THEOREM 5.1. *If $\eta_t(x) = 0$, then there exists a Toom graph $G = G(x,t)$, whose vertex set is included in $\{(z,s) : s \leq t, z \in \mathbf{Z}^2\}$ and which satisfies the following properties, for some sufficiently large $C > 0$:*

(1) *The number of possible graphs $G$ with $m$ edges is at most $C^m$.*
(2) *For a graph with $m$ edges there are at least $m/C$ vertices which are error sites.*
(3) *For any $r \geq 0$, if $\eta_t(y) = 0$ for $\|x - y\| \leq r$, then the number of edges of $G$ is at least $\max\{r/C, 3\}$.*

For the proof see [35].

In the classical application of Theorem 5.1, $b(x,t)$ are i.i.d. Bernoulli($p$). Then $P(\xi_t(x) = 1)$ converges to 1 as $p \to 1$, uniformly in $(x,t)$. Thus $\xi_t$ has an invariant measure with density close to 1. This also follows when $b(x,t)$ are not independent, but have uniformly bounded range of dependence in spacetime.

LEMMA 5.2. *Assume $u_1$ and $u_2$ are as in Lemma 4.5, with corresponding wedge $W = W_{u_1, u_2}$. Start a p-perturbation of $\mathcal{T}$ from $A_0 = W \cap \mathbf{Z}^2$. Consider the following two events:*

(1) $E_{x,t,M} = \{x \text{ is within distance } M \text{ of } A_t\}$.
(2) $F = \{\text{there is a } C \text{ so that, within the lattice ball of radius } t^2, (t - C\log t)z + W \subset A_t\}$.

*Then for any $\varepsilon > 0$, there exists an $M$ so that for any point $x \in tz + W$ the event $E_{x,t,M}$ happens with probability $1 - \varepsilon$ (uniformly in $x$ and $t$). Moreover, $P(F) = 1$.*

PROOF. It is convenient to translate the dynamics so that $W$ is fixed, that is, consider $A'_t = (A_t - tz) \cap W$. Also, rotate the lattice so that $W$ has its maximum at the origin and $u_1$ and $u_2$ are situated symmetrically with respect to the $y$-axis. It then follows from Corollary 4.3, Lemmas 4.4 and 4.5 that there are finite constants $C$ and $t_0$ (which again only depend on $\mathcal{T}$) so that the construction in the following paragraphs is possible.

Cut a finite neighborhood of the origin with a horizontal line $y = -C$, and let $t_0$ be the time the deterministic dynamics needs to fill $W$ again if sites above the cut are removed. Run the random system in multiples of time $t_0$, with the proviso that if any site at time $nt_0$ is 0 above the cut, then all sites above the cut are set to 0 immediately. Also make all the sites above the cut 0 if during the time interval $[(n-1)t_0 + 1, nt_0]$ a site within $C$ of the sites of the cut does not become occupied because of a bad coin flip, that is, although the deterministic dynamics would make it occupied the random one does not. The resulting set of occupied points at time $nt_0$ is called $A'_n$.

If an integer site $x \in W$ is not in $A'_n$, then either an integer site in $W$ strictly below the horizontal line through $x$ must be 0 in $A_{n-1}$ and a site on or above the horizontal line through $x$ must also be 0 in $A_{n-1}$, or else a site within distance $C$ of $x$ does not become occupied although the deterministic dynamics would make it occupied. In the latter case we call $x$ an *error site*. It is clear that error sites have finite range of dependence in spacetime and occur with probability $p'$ which is above a fixed power of $p$ and thus can be made arbitrarily close to 1.

By Theorem 5.1, uniformly for sites $x \in W$ and $n$,

$$P(x \text{ is at distance at least } M \text{ from } A'_n) \leq \sum_{m \geq M/C} C^m (p')^{m/C} \leq (Cp)^{M/C}.$$
(5.1)

Now the claim concerning the event (1) readily follows. To prove that $P(F) = 1$, note that we can choose $M = C \log n$ (this $C$ does depend on $p$), so that the probability in (5.1) is below $1/n^4$ and thus

$$P(\text{some } x \in W \cap B(0, n^2) \text{ is at distance at least } C \log n \text{ from } A'_n) \leq C/n^2.$$

From local regularity and Lemma 2.1, it now follows that

$$P((tz+W) \cap B(0,t^2) \not\subset A_{t+C\log t}) \leq C/t^2,$$

and Borel–Cantelli completes the proof. □

PROOF OF THEOREM 2 IN CASES 1 AND 2. Construct a convex set $\tilde{L}_n$ in the following manner:

Step 1. Start from $nL$, a multiple of the shape of the deterministic model.

Step 2. Every corner of $nL$ whose corresponding open edge of $K_{1/w}$ contains directions $u$ such that $u/w(u) \in \partial K'$ is "logarithmically rounded off" by introducing for each such $u$ an edge with normal $u$ and length $C \log n$ (where $C$ is some large constant whose value will become clear below), and ensuring that the resulting set is a convex subset of $nL$.

Step 3. Round off every corner of the set obtained in Step 2 to produce locally a translate of $W = W_{u_1,u_2}$ of Lemma 4.5.

If $R$ is the number of directions $u$ for which $u/w(u) \in \partial(\text{co}(K_{1/w}))$, then Step 3 has produced $R$ wedges $W$, which we label $W_0, \ldots, W_{R-1}$. Start with some large $\tilde{L}_M$ and couple (by using the same $\xi_{x,t}$) the resulting $R+1$ dynamics: one started from integer sites in each wedge, and the last one started from those inside $tL_M$.

If the percolation model introduced in the proof of Proposition 1.2 survives for all time on each edge of $\tilde{L}_t$, then we call the coupling *successful*. This means that the state of each site is *exactly the same* as the state of the same site in one of the wedges.

Fix an $x \in \tilde{L}_t$. By the FKG inequality and Lemma 4.7, the following three events simultaneously happen with probability at least $1 - \varepsilon$ if $M$ is large enough:

(1) The coupling is successful.
(2) $E_{x,t,M}$ happens for a suitable wedge $W_i$ guaranteed by (1).
(3) $F$ happens for every wedge $W_i$, $i = 0, \ldots, R-1$.

For an arbitrary initial set which fills space, we once again use the fact that $\tilde{L}_M$ is covered in finite time to get (S1) (2) (3) in Case 1, and (S1) (3) in Case 2.

It remains to show that the a.s. deviations from a corner in Case 2 are at least logarithmic. Pick a corner $\sigma \in L$ whose corresponding open edge of $\text{co}(K_{1/w})$ contains $u/w(u)$, for some direction $u$. Note that the boundaries of $nL$ and $nw(u)u + H_u^-$ intersect at exactly $t\sigma$. Finally, consider the infinite wedge $W$ defined by the corner and locate its vertex at the origin.

The number of sites in

$$S_k = W \cap ((-kw(u)u + H_u^-) \setminus (-(k-1)w(u)u + H_u^-))$$

RANDOM GROWTH MODELS 27is bounded above by $Ck$. Let $T_k$ be a geometric random variable with success probability $q_k = (1-p)^{|S_k|}$.

Start the dynamics from sites in $ML$, where $M$ is arbitrary. It is clear that

$$
\begin{aligned}
P(t\sigma \text{ is at distance at least } ck \text{ from } A_t) &\geq P(T_1 + \cdots + T_k \leq t) \\
&\geq P(T_k^{(1)} + \cdots + T_k^{(k)} \leq t),
\end{aligned}
\tag{5.2}
$$

where $T_k^{(i)}$ are i.i.d. copies of $T_k$. Now by Chebyshev

$$
P(T_k^{(1)} + \cdots + T_k^{(k)} \geq t) \leq \frac{k \operatorname{Var}(T_k)}{(t - kE(T_k))^2} = \frac{k}{q_k(tq_k - k)^2},
\tag{5.3}
$$

and it is easy to check that when $k = c\log t$ for a sufficiently small $c = c(p)$ the upper bound in (5.3) is $O(t^{-3/2})$. The desired result now follows from (5.2) and Borel–Cantelli. $\square$

We note that logarithmic a.s. fluctuations (S3) are optimal in Case 1 of Theorem 2 as well, as any of the sites in $(t+1)L \setminus tL$ can stay unoccupied for time $c\log t$ as a result of bad coin flips. This happens independently for each such site with probability $t^{-1/2}$ if $c = c(p)$ is small enough. Since the number of such sites is linear in $t$, a large deviation computation shows that the probability that this happens for at least one site is at least $1 - \exp(-c\sqrt{t})$ [19].

**6. An example.** In this section we present an example of a one-parameter family of random rules $\pi_p$, $p \in [0, 1]$, for which we can compute the half-space velocities explicitly. Every such example seems to be similarly based on the models introduced in [14] and [30]. Apart from the exactly stable cases, the example which follows therefore seems to be the *only* nontrivial instance of a random growth model with known shape.

The best way to think about this model is on the hexagonal lattice, but we will describe it so that it fits into our $\mathbf{Z}^2$ setup. The model's neighborhood $\mathcal{N}$ consists of seven sites, the von Neumann neighborhood with two added diagonal sites: $(1, -1)$ and $(-1, 1)$. Then $\pi_p(S) = 1$ when at least one of the following four conditions is satisfied: $(-1, 0) \in S$, $(0, 1) \in S$, $\{(1, -1), (0, -1)\} \subset S$, $\{(-1, 1), (1, 0)\} \subset S$. On all other nonempty sets $S$, $\pi_p(S) = p$. This is a $p$-perturbation of the additive model, but not a standard one. Nevertheless we will denote the half-space velocities by $w_p$ and the shapes by $L_p$. Supercriticality and local regularity are trivial.

Note that $\pi_p$ interpolates between two supercritical growth models. When $p = 1$, the CA is additive with neighborhood $\mathcal{N}$, which thus has $K_{1/w_1} = \mathcal{N}^*$ [which is $\operatorname{co}(\mathcal{N})$ rotated by $90°$], and $L_1 = \operatorname{co}(\mathcal{N})$. When $p = 0$, the vertices of $K_{1/w_0}$ are $(0, \pm 1)$, $(-1, 1)$, $(1, -1)$, $(1, 2)$ and $(-1, -2)$, while $L_0$ is a parallelogram with vertices at $(\pm 1, 0)$, $(-1/3, 2/3)$ and $(1/3, -2/3)$. What sets this model apart from a generic example is that certain initial sets make it *exactly solvable* for every $p$, in the sense that $P(x \in A_t)$ can be expressed as a Fredholm determinant of an explicitly known operator on $\ell^2$ [14]. These initial sets are four wedges, which together cover the plane: $W_1 = H_{e_2}^- \cap H_{-e_1}^-$, $W_2 = H_{e_2}^- \cap H_{e_d}^-$, $W_3 = H_{-e_2}^- \cap H_{e_1}^-$, $W_4 = H_{-e_2}^- \cap H_{-e_d}^-$, where $e_d = (1,1)/\sqrt{2}$.

Assume that $A_0 = W_1 \cap \mathbf{Z}^2$. The first observation is that $A_t = \{(x, y): y \le g_t(x)\}$, where $g_t: \mathbf{Z} \to [-\infty, t]$ is a nondecreasing function. This is easily proved by induction. As a consequence, whenever $x$ has its east (i.e., right) neighbor in $A_t$, it also has its southeast diagonal neighbor in $A_t$. On this initial set, the rule is therefore symmetric across the line $y = -x$. Now let, for every positive integer $n$, $h_t(n) = g_t(-t+n) - n$. Then $h_t(n) = -\infty$ for $n < 0$, $h_t(0)$ is a random walk which jumps by $+1$ (resp. stays put) with probability $p$ (resp. $1-p$), and for $n > 0$ $h_{t+1}(n)$ equals $h_t(n) + 1 = h_t(n-1)$ automatically whenever $h_t(n-1) > h_t(n)$. Finally, when $h_t(n-1) \le h_t(n)$, $h_{t+1}(n) = h_t(n) + 1$ with probability $p$ and otherwise $h_{t+1}(n) = h_t(n) + 1$. This establishes the equivalence. It follows that there exists a self-invertible function $\phi: [p, 1] \to [p, 1]$, so that $A_t/t$ converges to the region $L_{W_1} = \{(x, y) \in \mathbf{R}^2: y \le 1, -\phi(y) \le x$ for $p \le y$, and $-1 \le x$ for $y \le p\}$. The function $\phi$ has the following explicit form [14]:

$$\phi(y) = 1 - p - (1 - 2p)y + 2\sqrt{p(1-p)y(1-y)}.$$

The argument is similar when $A_0 = W_2 \cap \mathbf{Z}^2$. (In fact, that this case is equivalent to the above can be seen by mapping the model onto the hexagonal lattice, where it has fourfold symmetry.) If $g_t(x) = \sup\{-y + x \in A_t : y \in \mathbf{Z}\}$ (with $\sup \varnothing = -\infty$), then $h_t(n) = g_t(t-n) - n$ has the same evolution as the $h_t$ from the previous paragraph. It follows that this time $A_t/t$ converges to $L_{W_2} = \{(x, y) \in \mathbf{R}^2 : y \le 1, x \le \phi(y) - y$ for $p \le y$, and $x \le 1 - y$ for $y \le p\}$. The remaining two wedge shapes are obtained by symmetry: $L_{W_3} = -L_{W_1}$ and $L_{W_4} = -L_{W_2}$.

The proof of the next proposition is very similar to the proof of Corollary 1.1, and hence omitted. (See also [10].)

PROPOSITION 6.1. *Assume that a perturbation of a locally regular supercritical CA is given by $\pi$. Assume that its initial set is a wedge: $A_0 = W \cap \mathbf{Z}^2$, where $W = H_{u_1}^- \cap H_{u_2}^-$ and $u_1$ and $u_2$ form an angle in $(0, \pi)$. Then*

$$\frac{A_t}{t} \to \bigcap \{w_\pi(u)u + H_u^- : W \subset H_u^-\} = (K_{1/w_\pi} \cap W^*)^*,$$

*almost surely and in Hausdorff metric within any large ball of radius $C$.*

This proposition, together with Corollary 1.1, immediately implies that in the present example,
$$L_p = L_{W_1} \cap L_{W_2} \cap L_{W_3} \cap L_{W_4}$$
and therefore its top half comprises points $(x, y)$ which satisfy

$$-\phi(y) \leq x \leq \phi(y) - y, \quad \text{if } p \leq y \leq y_0,$$
$$-1 \leq x \leq 1 - y, \quad \text{if } 0 \leq y \leq p,$$

where

$$y_0 = y_0(p) = \begin{cases} 1, & p \geq 1/2, \\ \dfrac{2(1-p)}{3 - \sqrt{8p}}, & p < 1/2. \end{cases}$$

Hence the top half of the shape $L_p$ is convex and $\mathcal{C}^1$, but not strictly convex, for $p > 1/2$; is strictly convex and $\mathcal{C}^1$ only at $p = 1/2$; and for $p < 1/2$ is strictly convex with a corner at its highest point $(-y_0/2, y_0)$ above the $x$-axis. See Figure 5 for a plot.

We note that fluctuations from the limiting shape in every direction, except $\pm(-1/2, 1)$ when $p \leq 1/2$, can be obtained from [14]. For example, consider $\alpha \in (-1, (-y_0/2) \wedge (-p))$ and let $g_t(x) = \inf\{y \in \mathbf{Z} : (x, y) \in A_t\}$. (It is easy to see that sites in $A_t$ with a fixed $x$-coordinate always form an interval if they do so at $t = 0$.) Then $(g_t(\lfloor \alpha t \rfloor) - \phi(-\alpha)t)/t^{1/3}$ converges in distribution to a nondegenerate random variable. This follows from the fact that for such $\alpha$ the evolution of $A_t$ starting from a finite set and from $W_1 \cap \mathbf{Z}^2$ can be coupled so that the difference of respective $g_t(\lfloor \alpha t \rfloor)$ is stochastically bounded.

**7. Exact stability for box neighborhood TG CA.** The TG CA with box neighborhood of radius $\rho$ has $\rho(2\rho+1)$ supercritical thresholds $\theta$ [10], and the same number of corresponding $K_{1/w}$ which we label $K_1, \ldots, K_{\rho(2\rho+1)}$, and superimposed with different shades in Figure 6. Let $\mathcal{E} = \mathcal{E}_\rho = \bigcup_{\theta=1}^{\rho(2\rho+1)} \partial K_\theta$. At first this set appears to be of bewildering complexity (cf. the range 5 example on the right-hand side of Figure 6). However, perhaps the first feature revealed upon closer inspection is that $\mathcal{E}$ consists entirely of straight lines, called *K-lines*, which extend through the entire picture. [In fact, if we were to include the *critical* $K_{1/w}$ for $\theta = \rho(2\rho+1) + 1, \ldots, \rho(2\rho+1) + \rho$ all these lines would continue indefinitely.] There is one-to-one correspondence between the $K$-lines and the points of $\mathcal{N} \setminus \{0\}$, as we now explain.

For any of the $(2\rho+1)^2 - 1$ sites $x \in \mathcal{N} \setminus \{0\}$, start with a line through $x$ and $0$. Then rotate this line in the positive direction until it hits $0$ again. Call all such rotations $\ell_{x,\phi}$, $0 < \phi < 2\pi$. The set of all cardinalities $\Lambda(x, \phi) = |\mathcal{L}^-(\ell_{x,\phi})|$ is exactly the set of $\theta$ for which $x \in \ell_u$, for some $u$. Moreover, by Lemma 3.1, whenever $x$ is the only site in $\ell_{x,\phi} \cap \mathcal{N}$, this line determines a direction pointing toward the interior of an edge of $K_{\Lambda(x,\phi)}$, namely the direction of the normal to $\ell_{x,\phi}$ which points toward the origin. When the number of points in $\ell_{x,\phi} \cap \mathcal{N}$ to the left and to the right of $x$ are equal, the direction points to the interior of an edge in $K_\theta$ for $\theta = \Lambda(x, \phi) - |\{\text{points to the left of } x\}|$.

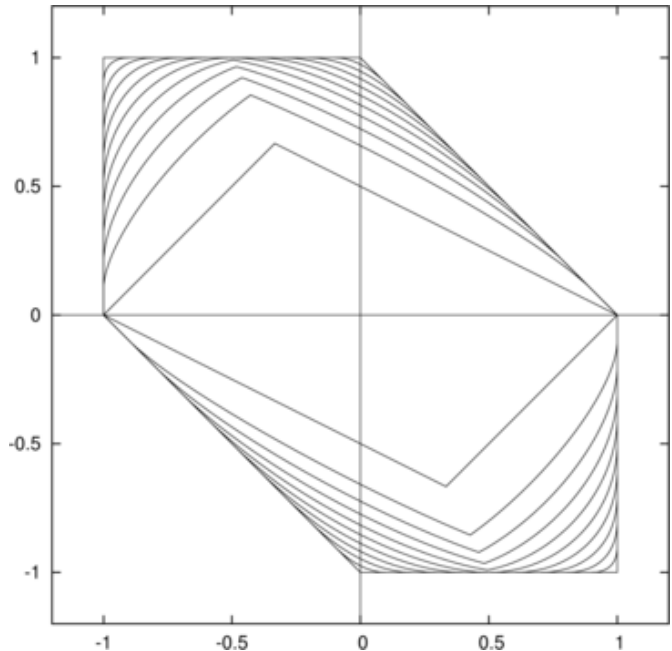

Fig. 5. *Shapes $L_p$ for $p = 0, 0.1, \ldots, 1$.*

Similarly, a line containing exactly two sites in $\mathcal{N} \setminus \{0\}$ determines a direction in which exactly two $\partial K_\theta$ meet (at a point which is a vertex of both, for one a convex vertex, for the other a concave one). A line containing exactly three points in $\mathcal{N} \setminus \{0\}$ determines a direction in which exactly three $\partial K_\theta$ meet (but this point is a vertex of only two of them). And so on. We are now ready to prove the next proposition, which in particular allows unambiguous reconstruction of all $K_\theta$ from $\mathcal{E}$.

PROPOSITION 7.1. *Two different $\partial K_\theta$ intersect only in a discrete set of points. Moreover, all finite tiles of $\mathcal{E}$ are triangles or quadrilaterals.*

PROOF. The first statement follows since for all but a discrete set of rotations $\phi$, $x$ is the only site in $\ell_{x,\phi} \cap \mathcal{N}$ and so the corresponding edge lies only in $K_{\Lambda(x,\phi)}$. Fix a $\theta_0$ and assume $\ell_{x,\phi} = \ell_u$ for some direction $u$. If $\ell_u$ contains two or more points in $\mathcal{N}$, then $\partial K_{\theta_0}$ does not intersect $\partial K_{\theta_0+1}$ if and only if $|\mathcal{L}^-(\ell_u)| = \theta_0$. Let $x$ be the rightmost point in $\ell_u \cap \mathcal{N}$ such that $\ell_{x,\phi} = \ell_u$ for some $\phi$. Decrease $\phi$ to the largest $\phi'$ such that $\ell_{x,\phi'} \cap \mathcal{N}$ contains at least two points. Clearly $\ell_{x,\phi'} \cap \mathcal{N}$ is a line $\ell_u$ for $\theta = \theta_0 + 1$ (and perhaps some larger $\theta$'s). A similar argument holds in the other direction, proving that, among two consecutive vertices of $\partial K_{\theta_0}$, at least one belongs to $\partial K_{\theta_0+1}$. This is clearly enough. □

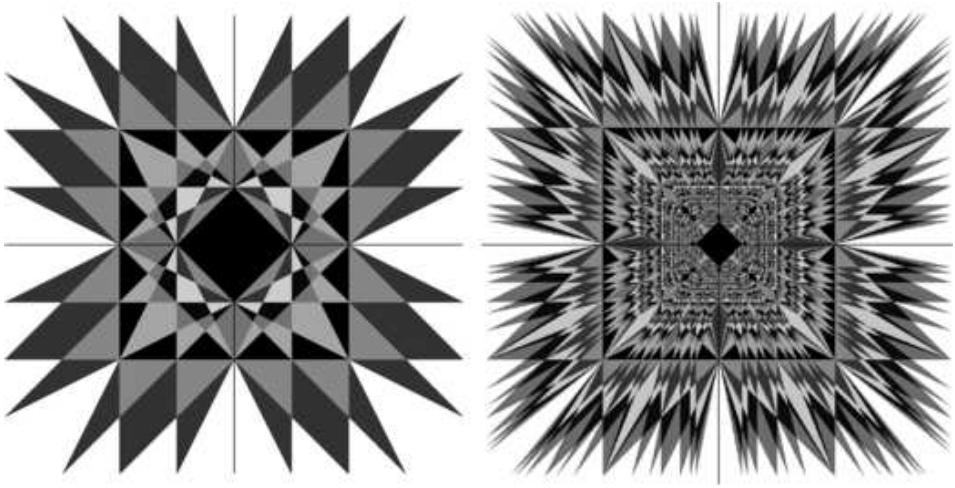

Fig. 6. *The* 10 *(resp.* 55*) supercritical* $K_{1/w}$*'s for range* 2 *(resp. range* 5*).*

For any $x$, let $\Lambda^*(x) = \inf_\phi |\Lambda(x,\phi)|$.

PROPOSITION 7.2. *$\partial K'$ includes a line segment if and only if $\theta = \Lambda^*(x)$, for some $x$. Moreover, for every $x$, $\Lambda^*(x) = |\Lambda(x,\phi)|$ for some $\phi$ such that $\ell_{x,\phi}$ exits $\tilde{\mathcal{N}} = B_\infty(0,\rho) \subset \mathbf{R}^2$ through the neighboring (rather than opposite) sides. Finally, each $K$-line includes exactly one line segment on $\partial K'$, for some $\theta$.*

PROOF. Pick an $x$ and $\phi_0$ so that $x$ is the only site in $\ell_{x,\phi_0} \cap \mathcal{N}$, and let $\theta = \Lambda(x,\phi_0)$. For this $\theta$ and normal $u$ to $\ell_{x,\phi_0}$, $\ell_u = \ell_{x,\phi_0}$. For the first assertion, it suffices to show that if $\theta > \Lambda^*(x)$, then $u/w(u)$ cannot lie on $\partial(\mathrm{co}(K_\theta))$. But if it would lie there, all of $K_{1/w}$ would have to lie on one side of the line of $K_{1/w}$ determined by $x$ and $\phi_0$ (see Lemma 3.1), meaning that $\Lambda(x,\phi) \geq \theta$ for every $\phi$, a contradiction.

For the second assertion, assume the said line exits the left- and right-hand sides of the square. We can also assume that $x$ lies either strictly inside the third quadrant, or on the negative $y$-axis. In both cases, rotate the line around $x$ to angle $\phi'$ in the negative direction just past the southeast corner of the square; this results in $\Lambda(x, \phi') \leq \Lambda(x, \phi)$ (with equality in the second case).

For the final claim, assume that $x$ is as in the above paragraph and that the line $\ell_{x,\phi}$ produces $\Lambda^*(x)$. Now rotate it in the negative direction to the smallest angle $\phi'$ for which the number of points on the left of $x$ on $\ell_{x,\phi'} \cap \mathcal{N}$ is larger than the number of points on the right. If $a$ and $b$ are the lengths of the line segments on $\ell_{x,\phi'}$ from the left and bottom edges of $\tilde{\mathcal{N}}$, respectively, then $a > b$. This remains the case for any $\phi'' < \phi'$ and thus further rotations in the negative direction only lose more points. An analogous argument works for positive rotations. □

Therefore, each $K$-line contributes exactly one edge on exactly one $\partial K'$. What produces a prevalence of exactly stable cases are those $\partial K'$ which have more than their share of edges. For example, it is immediate by symmetry that the number of edges on any $\partial K'$ is either 0 or at least 4. The number of $\theta$ with lack of exact stability hence does not exceed $((2\rho+1)^2 - 1)/4 = \rho^2 + \rho$, and therefore the number of exactly stable ones is at least $\rho^2$. This argument is quite simple, yet it fails to produce any way to identify a single exactly stable case. Our next result remedies this somewhat. However, we do not have an algorithm which lists more than $O(\rho)$ cases of either type.

PROPOSITION 7.3. *All $\theta \geq 2\rho^2 + 1$ are exactly stable cases. On the other hand, all $\theta \leq \rho$ and $\theta = \rho + 1 + i(2\rho + 1)$, $i = 0, \ldots, \rho - 1$, are not exactly stable.*

PROOF. For the second assertion, consider first $\theta \leq \rho$ and take $x = (-\rho + \theta - 1, -\rho)$. For any angle $\phi$, $\Lambda(x, \phi) \geq \theta$. This immediately implies that the edges in $K_{1/w}$ adjacent to the positive $y$-axis lie in $\partial K'$. If $\theta = \rho + 1 + i(2\rho + 1)$, $i = 0, \ldots, \rho$, there is instead a single edge, perpendicular to the $y$-axis.

For the first assertion, it is enough to prove that all first-quadrant boundary edges of such $K_\theta$ have slopes in $[1, \infty]$. Equivalently, take $x \in \mathcal{N}$ in the third quadrant strictly above the line through the origin with slope 1, and a line $\ell$ through $x$ with normal a negative rotation of $e_2$ by angle at most $\pi/4$ and such that $x$ is the only point in $\ell \cap \mathcal{N}$. Then $|\mathcal{L}^-(\ell)| < \theta$. This is certainly true when $u$ is close to vertical, and further rotations can only decrease $|\mathcal{L}^-(\ell)|$. □

Note that the proof of the above proposition shows that the $\theta = 2, \ldots, \rho$ thresholds have at least eight edges in $\partial K'$, which improves the lower bound on the number of exactly stable cases to $\rho^2 + \rho - 1$. This is a good lower bound for small $\rho$, although the exact enumerations are taxing.

PROOF OF THEOREM 3. By Proposition 7.2, we can reformulate the problem as follows. Consider all integer points $(a, b)$, $1 \leq a, b \leq \rho$. Take a line $\ell$ through $(a, b)$ which intersects the positive halves of both axes. For each such line, let $\Lambda(\ell)$ be the number of integer points in the closed triangle $T(\ell) \subset \mathbf{R}^2$ bounded by $\ell$ and the positive halves of the axes. Finally, let $\ell^* = \ell^*(a, b)$ be a line which minimizes $T(\ell)$. We need to find an upper bound for the number of different $T(\ell)$ over all $(a, b)$. In the following computations, $O(1)$ refers to a term of arbitrary sign whose absolute value can be bounded by a constant independent of $\rho$.

Let $\ell_x^*$ and $\ell_y^*$ be the line segments on $\ell^*$ from $(a, b)$ to the $x$- and $y$-axes, respectively. The first observation is that $\ell^*$ can be chosen so that the lengths of $\ell_x^*$ and $\ell_y^*$ differ by $O(1)$. (In fact, these lengths can be made arbitrarily close, although not necessarily equal.) In particular, the area of $T(\ell)$ is $2ab+O(1)$. We will assume, without loss of generality, that the length of $\ell_y^*$ is not larger than the length of $\ell_x^*$.

Now find an integer point $(0, y_0)$ on the $y$-axis immediately below where $\ell^*$ intersects the $y$-axis. Reflect $(0, y_0)$ through $(a, b)$ to get a point $(x_1, y_1)$ within $O(1)$ of the intersection of $\ell^*$ with the $x$-axis. [Note that $(x_1, y_1)$ lies within the closed first quadrant.] Also, let $(x_0, 0)$ be an integer point on the $x$-axis immediately to the left of where $\ell^*$ intersects it. Form the closed polygon $\Pi \subset \mathbf{R}^2$ by connecting $(a, b) \to (0, y_0) \to (0, 0) \to (x_0, 0) \to (x_1, y_1) \to (a, b)$.

The area of $\Pi$ is $2ab + O(1)$. Most importantly, the number of integer sites in $\Pi$ is $\Lambda(\ell^*) + n^*/2 + O(1)$, where $n^*$ is the number of integers in $\partial\Pi$ which are not on the axes. This follows since one loses about half of these by a small rotation of the line from $(0, y_0)$ to $(x_1, y_1)$ around $(a, b)$. Therefore,

$$T(\ell^*) = |\{\text{integer points in the interior of } \Pi\}| + \frac{n^*}{2} + 2a + 2b + O(1),$$

and by Pick's theorem [1],

$$\text{area of } \Pi = |\{\text{integer points in the interior of } \Pi\}| + \frac{n^*}{2} + a + b + O(1)$$
$$= T(\ell^*) - a - b + O(1).$$

It follows that

$$T(\ell^*) = 2ab + a + b + O(1) = \tfrac{1}{2}(2a+1)(2b+1) + O(1).$$

Let $M_N$ be the number of different products $mn$ of integers $m, n \in [1, N]$. It follows that

$$|\{T(\ell^*) : 0 \le a, b \le \rho\}| \le C \cdot M_{(2\rho+1)} \le C \cdot \frac{\rho^2}{\log^h \rho},$$

by the Hall–Tenenbaum sharpening of a theorem of Erdős ([17], Theorem 23). This ends the proof of the upper bound. The lower bound follows because the lower bound in [17], Theorem 23 is obtained using odd integers only. $\square$

**8. Final remarks.** In this section we mention several results which are related to the main topics of the paper, sketch their proofs, and also complete the proof of Theorem 1 by removing $\log^2 t$ from the large deviation estimate.

REMARK 1 (Continuous time). A standard continuous-time growth model $\tilde{A}_t$ is obtained by adjoining every site $x$ at an independent rate 1 exponential time after the time $\tau_x$ at which $x \in \mathcal{T}(\tilde{A}_{\tau_x})$. This process can be constructed in the standard way by attaching a Poisson process $\tilde{\xi}_x$ to every $x$. Theorem 1 is still valid in this case. The a priori large deviation bound however has $\log^2 t$ replaced by $\log^4 t$. We now sketch the proof.

Observe $\tilde{A}_t$ in discrete time units $t = 1, 2, \ldots$. Change it to $\tilde{A}'_t$ by making sure that between each time $t$ and $t+1$ no site at distance more than $C \log t$ from $\tilde{A}'_t$ gets occupied. As is easy to see by comparison with the continuous time additive dynamics having the same neighborhood,

$$P(\tilde{A}_t \neq \tilde{A}'_t \text{ within a lattice ball of radius } t^2) < t^{-3},$$

when $C$ is large enough. Now continue the proof with $\tilde{A}'_t$, which of course is a discrete-time monotone Markov process. Lemma 2.2 must be used $C \log t$ successive times to obtain the analog of estimate (2.2), and here is where the larger power of log originates. From this point the proof proceeds on familiar grounds, yielding existence of the asymptotic speeds $\tilde{w}$, while the continuous-time version of Corollary 1.1 establishes existence of the shape $\tilde{L}$.

REMARK 2 (Approximating half-space velocities). Perhaps the most convenient method for simulating a growth CA started from a half-space is to use a strip with tilted periodic boundary. Take a vector $u$ at angle $\phi \in [0, \pi/4]$ to $e_2$. (By rotations and reflections it is clearly enough to consider these.) Then take a large $L$, and restrict the growth to the strip $H_M = [0, M-1] \times \mathbf{Z}$. Let $\kappa = \tan \phi$. Given any configuration of occupied sites inside $H_M$, extend it to $\bar{H}_M = [-M, 2M-1] \times \mathbf{Z}$, by identifying the state of $(x, y)$ with that of $(x - M, \lfloor (y - \kappa M) \rfloor)$ if $x \geq M$ and with that of $(x + M, \lfloor (y + \kappa M) \rfloor)$ if $x < 0$.

Start from $A_0 = \{(x, y) \in H_M : y \leq \kappa x\}$. Let the dynamics update sites in $H_M$ with the specified boundary conditions until some large time $t$ when the interface apparently equilibrates. At this point, the average height above all points in $[0, M-1]$, multiplied by $\cos \phi / t$, is a good approximation to $w_\pi(u)$.

Note that this is a *much* more efficient technique for computing the shape $L_p$ than merely running the dynamics from a finite seed and observing the resulting blob. In particular, smoothness of $L_p$ is impossible to discern that way. The method outlined above, on the other hand, uses averaging to greatly reduce transversal fluctuations on the interface. The theoretical underpinning is partly given in our last theorem.

For a fixed $M$, and any $x \in [0, M-1]$, let

$$h_t^1(x) = \max\{y : (x, \lfloor y - \kappa x \rfloor) \in A_t\}, \qquad h_t^1 = \max\{h_t(x) : x \in [0, M-1]\},$$
$$h_t^0(x) = \min\{y : (x, \lfloor y - \kappa x \rfloor) \notin A_t\}, \qquad h_t^0 = \min\{h_t(x) : x \in [0, M-1]\}.$$

THEOREM 4. *Fix an $\varepsilon > 0$. If $M$ is large enough, then with probability 1*

$$\frac{w_\pi(u) - \varepsilon}{\cos \phi} \leq \liminf_{t \to \infty} \frac{h_t^0}{t} \leq \limsup_{t \to \infty} \frac{h_t^1}{t} \leq \frac{w_\pi(u) + \varepsilon}{\cos \phi},$$

*uniformly in $u$.*

In fact, it is easy to show by subadditivity that $h_t^0/t$ and $h_t^1/t$ both converge a.s. to the same number as $t \to \infty$.

PROOF OF THEOREM 4. We start by proving the lower bound. Note that the boundary effects spread with finite speed, as $\mathcal{N}$ is finite. Thus, until the time $t_0 = cM$, the occupied sites on any vertical line through $x \in [0, M-1]$ are above those started from an infinite tilted half-plane through $(x, \lfloor \kappa x \rfloor)$ or through $(x, \lfloor \kappa x \rfloor - 1)$. (We have to allow for the second possibility because there may not be a perfect match at the boundaries.) By the weaker form of Theorem 1, the probability that the lowest unoccupied site above a fixed $x$ is below $\kappa x + (w_\pi(u)(1-\varepsilon)/\cos \phi)t_0$ is at most $\exp(-ct_0/\log^2 t_0)$. For $n \geq 0$, define this translation of $A_0$:

$$B_n = \left\{(x, y) : y \leq \kappa x + \frac{w_\pi(u) - \varepsilon}{\cos \phi} n t_0\right\}$$

and the event

$$E = \{B_1 \subset A_{t_0}\}.$$

It follows that

$$P(E^c) \leq M \exp(-ct_0/\log^2 t_0) < \varepsilon,$$

if $M$ is large enough.

Now run the dynamics until time $t_0$. If $E$ happens, restart the dynamics from the set $B_1$; otherwise restart the dynamics from $B_0 = A_0$. Then repeat from the possibly translated $A_0$. Let $U_n$ be the largest $u$ for which $B_u \subset A_{nt_0}$. We have just proved that $U_n$ dominates an $n$-step random walk which at each step increases by 1 with probability $1 - \varepsilon$ and stays put with probability $\varepsilon$. Therefore $U_n \geq (1-2\varepsilon)n$ with probability at least $1-\exp(-cn)$. This ends the proof of the lower bound (as uniformity in $u$ follows because the constant $c_\varepsilon$ in the weaker form of Theorem 1 does not depend on $u$).

The upper bound is proved similarly, except for the fact that we need an upper bound on the extent to which $A_t$ can propagate in $t_0$ steps. The trivial bound $Ct_0$, where $C$ is the diameter of $\mathcal{N}$, suffices. The comparison random walk now increases by 1 with probability at least $1-\varepsilon$ and by $C$ with probability $\varepsilon$, and so its speed is $1+C\varepsilon$. □

Theorem 4 remains valid for continuous-time growth as well, with a similar proof, except for two significant differences. The first is that the coupling of finite and infinite systems only holds up to time $cL$ with probability exponentially close to 1 in $L$ [15]. The second is that the trivial upper bound at end of the proof is not available, so the forward jump of the comparison random walk is arbitrarily long, with exponential tail probabilities. It is clear that the proof is still valid under these conditions.

PROOF OF THEOREM 1 (CONCLUDED). We choose $u$, $M$ and $t_0$ as in the proof of Theorem 4. Recall that $t_0/M$ is small and so during a time interval of length $t_0$ sites at a distance larger than $M$ cannot interact. Define

$$B_{n,i} = \left\{(x,y) : iM \leq x \leq (i+1)M - 1 \text{ and } y \leq \kappa x + \frac{w_\pi(u) - \varepsilon}{\cos\phi} nt_0\right\}$$

and the event

$$E = \{B_{1,0} \subset A_{t_0}\}.$$

As before, $P(E^c) < \varepsilon$ if $M$ is large enough. Now we define $\tilde\eta_n(i)$ to be the maximal $k$ such that every site $(x,y)$ with $x \in [iM, (i+1)M - 1]$ and $y \leq \frac{w_\pi(u)-\varepsilon}{\cos\phi} kt_0 + \kappa x$ is in $A_{nt_0}$. We can couple $\tilde\eta_n$ and the slow version of $\eta'_n$ from Lemma 3.3 (with $p' = 1 - \varepsilon$), so that $\tilde\eta_n(i) \geq \eta'_n(i)$ for every $i$ and $n$. (The Bernoulli random variables $b$ are probabilities of suitable translates of the event $E$.) It follows that, for every $\delta > 0$ we can find a small enough $\varepsilon$ (which then dictates a large enough $M$), so that

$$P\left((0,y) \in A_t \text{ for all } y \leq (1-\delta)\frac{w_\pi(u) - \varepsilon}{\cos\phi} t\right) \geq 1 - e^{-ct}.$$

This proves the exponential bound on the probability that $A_t/t$ lags significantly behind $w_\pi(u)$. We will now show that it cannot progress significantly faster than $w_\pi(u)$ either. To this end, we redefine

$$B_{n,i} = \left\{(x,y) : iM \leq x \leq (i+1)M - 1 \text{ and } y \leq \kappa x + \frac{w_\pi(u) + \varepsilon}{\cos\phi} nt_0\right\}$$

and

$$E = \{A_{t_0} \subset B_{1,0}\},$$

so that again $P(E^c) < \varepsilon$. In case $E$ fails, the occupied sites above the interval $[0, M-1]$ cannot progress by more than the diameter of $\mathcal{N}$. Furthermore, now $\tilde{\eta}_n(i)$ is the maximal $k$ such that some site $(x,y)$ with $x \in [iM, (i+1)M-1]$ and $y \leq Rk + \frac{w_\pi(u)+\varepsilon}{\cos\phi}nt_0 + \kappa x$ is in $A_{nt_0}$. Here $R$ is a suitably large multiple of the diameter of $\mathcal{N}$, which ensures that the fast version of $\eta'_n$ in Lemma 3.3 (now with $p' = \varepsilon$) dominates $\tilde{\eta}_n$. Thus Lemma 3.3(3) completes the proof. $\square$

REMARK 3 (Continuity of $K_{1/w_p}$). Assume a standard $p$-perturbation of $\mathcal{T}$. As $p$ changes from 0 to 1, $K_{1/w_p}$ varies continuously. To see this, assume that $p'$ is close to $p$, $p' < p$ and couple the systems with the two probabilities in the obvious way. To show that $w_{p'}(u)$ is close to $w_p(u)$ (uniformly in $u$), one needs to look at the proof of the lower bound in Theorem 4. Between 0 and $t_0$, the occupied sets in the two systems will not differ at all with probability $(1-(p-p'))^{Ct_0}$, which, as $t_0$ is constant (albeit dependent on $\varepsilon$), can be made larger than $1-\varepsilon$ if $p-p'$ is small enough. Once this observation is made, it is only necessary to follow the rest of the lower bound proof with $p$ replaced by $p'$.

Note that this continuity alone demonstrates that $L_p$ has corners (i.e., is not differentiable) in every case not quasi-additive when $p$ is close enough to 1 (although these corners may not move at the same speed, or even in the same direction, as the corresponding corners of $L_1$).

REMARK 4 (Shapes for small $p$). Again, assume a standard $p$-perturbation of $\mathcal{T}$. What happens as $p \to 0$? Certainly $L_p$ shrinks, and in fact

$$\frac{1}{p}L_p \to \tilde{L},$$

the limit shape of the continuous-time growth model $\tilde{A}_t$. To see this, let $A'_t$ be $A_{\lfloor t/p \rfloor}$. After a site sees a sufficient configuration in $A'_t$ (resp. $\tilde{A}_t$), it becomes occupied in a time distributed as $T_g$ (resp. $T_e$). It is easy to see that for small $p$, distributionally, $T_g \geq T_e(1-p)$. Rescaling of continuous time immediately gives $\tilde{L}/(1-p) \supset L_p/p$ for small $p$.

For the opposite direction, observe that $T_e + p \geq T_g$, in distribution. The lower bound part of the proof of Theorem 4 for continuous time now shows that $\tilde{A}_{t_0} \supset B_1$ with probability $1-\varepsilon$. With the same probability, then, $A'_t \supset B_1$ at time $t = t_0 + CpLt_0 < (1+\varepsilon)t_0$ if $p$ is small enough. (The added term is simply $p$ times the number of sites in $B_1 \setminus A_0$.) This easily completes the proof.

In closing, let us pose two challenging conjectures based on experiment.

CONJECTURE 8.1. *If $\pi$ is a standard p-perturbation of a locally regular and supercritical CA with convex $K_{1/w}$, then $K_{1/w_p}$ is strictly convex.*

This conjecture fails for nonstandard perturbations, as seen from the example discussed in Section 6. In fact, it fails for nonstandard perturbations even if we restrict to $p$ very close to 1. A range 2 box counterexample is obtained by $\pi(S)$ which is 1 when $|S| \geq 9$ and $p$ when $|S| = 8$. A glance at Figure 6, together with results from Section 3, confirms that $K_{1/w_\pi}$ cannot be convex for any $p < 1$.

Figure 7 illustrates the application of Theorem 4 to two examples. The left frame depicts $K_{1/w_p}$ for the Moore TG CA with $\theta = 3$ and $p = 1, 0.9, \ldots, 0.4$, while the right frame does the same for the range 2 box TG CA with $\theta = 9$ and $p = 1, 0.975, \ldots, 0.5$. As guaranteed by Theorem 2, $\text{co}(K_{1/w_p}) \subset \text{co}(K_{1/w_1})$ for $p$ close enough to 1. What is more, the angles at the corners of $K_{1/w_p}$ approach the angles of $K_{1/w_1}$ as $p \to 1$. (This can in fact be proved by methods of the present paper.) As $p$ decreases, $\partial K_{1/w_p}$ separates from $\partial K_{1/w_1}$, but $K_{1/w_p}$ remains nonconvex. (Whether it may become $\mathcal{C}^1$ before the boundaries separate is not clear.) Upon further decrease in $p$ one observes concavities gradually filling in, until $K_{1/w_p}$ becomes convex. Such observations, as well as early belief in the asymptotic isotropy of Eden's continuous-time random growth model [7], suggest our final conjecture.

CONJECTURE 8.2. *If $p$ is small enough, the standard p-perturbation of $\mathcal{T}$ has strictly convex and smooth $K_{1/w_p}$. The continuous-time version $K_{1/\tilde{w}}$ is also strictly convex and smooth.*

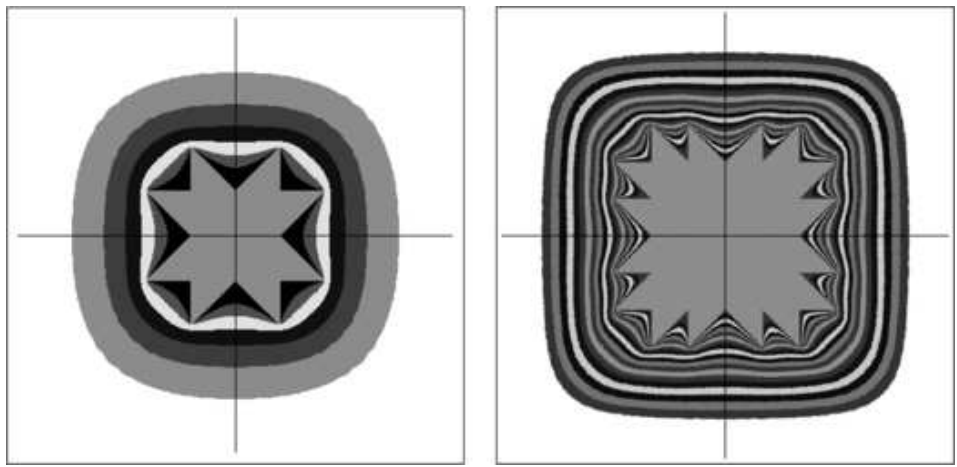

FIG. 7. *Two examples of TG $K_{1/w_p}$'s.*

**Acknowledgments.** We thank Gérald Tenenbaum for kindly sharing his expertise on intricacies of the Linnik–Vinogradov–Erdős problem with us. Thanks also go to Timo Seppäläinen for pointing out an error in an early version.

Mathematics Department  
University of California  
Davis, California 95616  
USA  
e-mail: gravner@math.ucdavis.edu

Department of Mathematics  
University of Wisconsin  
Madison, Wisconsin 53706  
USA  
e-mail: griffeat@math.wisc.edu